\documentclass[11pt,a4paper]{amsart}

\usepackage[pdftex, 
	a4paper,
	unicode,
	linktocpage,
	urlbordercolor={1 0.5 0.125},
	filebordercolor={1 0.5 0.125},       
	linkbordercolor={0.25 1 0.25},
	citebordercolor={0.25 1 0.25},
	pdfusetitle,
	pagebackref=false,
	pdfstartview=FitH,
	pdfborderstyle={/S/U/W 1}, 
	breaklinks]{hyperref}

\usepackage[T1]{fontenc}
\usepackage{lmodern}
\usepackage[utf8]{inputenc}

\usepackage{graphicx}
\usepackage{amssymb}

\usepackage{verbatim} 

\newtheorem*{TheoremA}{Theorem A}
\newtheorem*{TheoremB}{Theorem B}
\newtheorem*{TheoremC}{Theorem C}
\newtheorem{theorem}{Theorem}[section]
\newtheorem{proposition}[theorem]{Proposition}
\newtheorem{lemma}[theorem]{Lemma}
\newtheorem{corollary}[theorem]{Corollary}

\theoremstyle{definition}
\newtheorem*{definition}{Definition}
\newtheorem{question}{Question}

\theoremstyle{remark}
\newtheorem{remark}[theorem]{Remark}

\newcommand{\ie}{i\textup.e\textup.{}}

\newcommand{\Int}[2]{\mathfrak{Int}_{#1}(#2)}
\newcommand{\id}{\operatorname{id}}
\newcommand{\fix}{\operatorname{Fix}}
\newcommand{\supp}{\operatorname{Supp}}

\newcommand{\Th}{\operatorname{Th}}
\newcommand{\Mod}{\operatorname{Mod}}
\newcommand{\arity}{\operatorname{ar}}
\newcommand{\V}{\mathcal V}
\newcommand{\T}{\mathcal T}
\newcommand{\F}{\mathcal F}

\title[Arithmetic in groups of permutations of an interval]
{Interpretation of the Arithmetic\\
in certain groups of piecewise affine\\
permutations of an interval}

\author[T.~Alt\i nel \and A.~Muranov]{Tuna Alt\i nel \and Alexey Muranov}

\date{\today}


\address{Universit\'e de Lyon;
Universit\'e Lyon 1;
Institut Camille Jordan CNRS UMR 5208;
43 boulevard du 11 novembre 1918;
F--69622 Villeurbanne Cedex, France}


\thanks{The second author was supported first by
Chateaubriand Fellowship and later by 
\emph{FP6 Marie Curie Research Training Network
in Model Theory and its Applications\/}
funded by the European Commission under contract number
\texttt{MRTN-CT-2004-512234} (MODNET)}

\begin{document}



\begin{abstract}
The Arithmetic is interpreted in all the groups of Richard Thompson
and Graham Higman, as well as in other groups of piecewise affine
permutations of an interval which generalize the groups of Thompson
and Higman.
In particular, the elementary theories of all these groups are
undecidable.
Moreover, Thompson's group $F$ and some of its generalizations
interpret the Arithmetic without parameters.
\end{abstract}

\maketitle

\tableofcontents



\section{Introduction}
\label{section:introduction}

Valery Bardakov and Vladimir Tolstykh
\cite{BardakovTolstykh:2007:iaTgF} 
showed recently that Richard Thompson's group $F$
interprets the Arithmetic.
In other words, $F$ interprets the structure $(\mathbb N,+,\times)$
by first-order formulae with parameters.
In this work we generalize their result in two directions.
On the one hand, in Sections \ref{section:twolemmas} and
\ref{section:copiesZwrZ},
we generalize the approach of Bardakov and Tolstykh to
the construction of the Arithmetic from $F$, and show
that it works for all the groups defined by
Melanie Stein in \cite{Stein:1992:gplh},
which include all the three groups of Thompson and all the groups of
Graham Higman \cite{Higman:1974:fpisg}.
On the other hand, in Section \ref{section:interpretationarithmetic}
we show that the group $F$ and some of its generalizations
interpret the Arithmetic \emph{without parameters}.
(The difference between \emph{with\/} et \emph{without\/} parameters
shall be explained in Section~\ref{subsection:generalities.modeltheory}.)

The elementary theory of the Arithmetic $(\mathbb N,+,\times)$
is famous for its complexity since Gödel's incompleteness theorems
\cite{Godel:2006:fuSPMvS1-ger}.
If a structure of finite signature interprets the Arithmetic,
then the elementary theory of this structure is
\emph{hereditarily undecidable}.
(A theory of finite signature is called
\emph{hereditarily undecidable\/} if
every its subtheory of the same signature is undecidable,
see \cite[\S3]{Tarski:1971:gmpu}.)
Indeed, it is well known from the work Andrzej Mostowski, Raphael Robinson,
and Alfred Tarski
\cite{MostowskiRobTar:1971:ueua,MostowskiTarski:1949:uaitr,
Tarski:1971:gmpu} that the elementary theory of the Arithmetic is
hereditarily undecidable.
It is also well known to specialists that if a
structure $N$ of finite signature interprets with parameters another
structure $M$ of finite signature whose elementary theory is
hereditarily undecidable,
then the elementary theory of $N$ is hereditarily undecidable as well.%
\footnote{In \cite{BardakovTolstykh:2007:iaTgF} the authors
state this fact with a reference to \cite{Ershov:1980:prkm-rus}.
We show this fact as Lemma \ref{lemma:19.(6.2)}.}
Thus Bardakov and Tolstykh showed that
the elementary theory of $F$ is hereditarily undecidable,
and hence Question 4.16 by Mark Sapir in
\cite{anonymous:2004:Tg40y} is partially resolved.
The same argument shows that
the elementary theories of all the groups that we study
in this paper are hereditarily undecidable.

For the reader's convenience, we present
in Section \ref{section:undecidability} our version of a proof,
also based on a result of Mostowski, Robinson, and Tarski
\cite[Theorem 9]{MostowskiRobTar:1971:ueua},
that if a structure $S$ of finite signature
interprets the Arithmetic with parameters,
then the elementary theory of $S$ is hereditarily undecidable.

The groups that are the object of our study appear naturally
as generalizations of the tree groups defined by Thompson in 1965
and customarily denoted $F$, $T$, and $V$.%
\footnote{Other letters were also used to denote these groups
(see \cite{CannonFloydParry:1996:inRTg}).
It is rather common, for example, to denote the group $V$ by~$G$.}
Thompson's groups are presented in detail in
\cite{BelkBrown:2005:fdeTgF,CannonFloydParry:1996:inRTg}.
All the three are infinite finitely presented.
The group $V$ was the first known example of a finitely presented
infinite simple group.
The group $T$ is also simple.
The group $F$ embeds in $T$, and $T$ embeds in $V$.
The group $V$ was generalized by Higman \cite{Higman:1974:fpisg}
into a family of finitely presented groups
$G_{n,r}$, $n=2,3,4,\dotsc$, $r=1,2,3,\dotsc$,
where $G_{2,1}\cong V$.
Higman's group $G_{n,r}$ is simple if $n$ is even;
if $n$ is odd, the derived subgroup $[G_{n,r},G_{n,r}]$ is
simple of index $2$ in $G_{n,r}$.
Kenneth Brown \cite[Section 4]{Brown:1987:fpg}
similarly generalized the groups $F$ and~$T$.

Thompson's groups have representations by piecewise affine permutations
of an interval, where the group $F$ is represented by
homeomorphisms with respect to the usual topology,
and $T$ is represented by
homeomorphisms with respect to the topology of circle.
Stein \cite{Stein:1992:gplh} studied three families of groups of
such permutations which generalize respectively the
three groups of Thompson.
In order to state our main results, we shall review here
the definitions of these families.

Let $r$ be a positive real number and $\Lambda$ a subgroup of
the multiplicative group $\mathbb R^*_+$ of positive reals.
Let $A$ be an additive subgroup of $\mathbb R$ containing $r$
and invariant under the action of $\Lambda$ by multiplications.
Then define $\V(r,\Lambda,A)$ to be the group of all the
bijections $x\colon [0,r)\to[0,r)$ that satisfy the following conditions:
\begin{enumerate}
\item
	$x$ is piecewise affine with finitely many
	cuts and singularities;
\item
	$x$ is right-continuous at every point (in the usual sense);
\item
	the slope of each affine part of $x$ is in $\Lambda$;
\item
	all cut and singular points of $x$, as well as their images,
	are in~$A$.
\end{enumerate}
The family of groups $\V(r,\Lambda,A)$ contains all the groups of Higman:
for every $n=2,3,\dotsc$ and every $r=1,2,\dotsc$,
$$
G_{n,r}\cong\V(r,\langle n\rangle,\mathbb Z[{\textstyle\frac{1}{n}}]).
$$

Define subgroups $\F(r,\Lambda,A)$ and $\T(r,\Lambda,A)$
of $\V(r,\Lambda,A)$ as follows:
\begin{itemize}
\item
	$\F(r,\Lambda,A)$ is the subgroup of all the elements of
	$\V(r,\Lambda,A)$ continuous with respect to the usual
	topology of $[0,r)$,
\item
	$\T(r,\Lambda,A)$ is the subgroup of all the elements of
	$\V(r,\Lambda,A)$ continuous with respect to the topology
	of circle on $[0,r)$
\end{itemize}
(where the topology of circle on $[0,r)$ is the topology
induced by the natural identification of $[0,r)$ with
the topological quotient $[0,r]/\{0,r\}$).
The groups $F$, $T$, and $V$ of Thompson are isomorphic to
$\F(1,\langle 2\rangle,\mathbb Z[1/2])$,
$\T(1,\langle 2\rangle,\linebreak[0]\mathbb Z[1/2])$,
and $\V(1,\langle 2\rangle,\mathbb Z[1/2])$, respectively.
Groups of the form $\F(r,\Lambda,A)$ were studied already by
Robert Bieri and Ralph Strebel in 
\cite{BieriStrebel:pp1985:gPLhrl} (unpublished).

For the rest we shall always assume that $\Lambda\ne\{1\}$.

\begin{TheoremA}
If\/ $G$ is a subgroup of\/ $\V(r,\mathbb R^*_+,\mathbb R)$
such that
$$
G\cap\F(r,\mathbb R^*_+,\mathbb R)=\F(r,\Lambda,A),
$$
then\/ $G$ interprets the Arithmetic\/
$(\mathbb N,+,\times)$ with parameters\textup.
\end{TheoremA}

\begin{TheoremB}
If\/ $\Lambda$ is cyclic\textup,
then\/ $\F(r,\Lambda,A)$ interprets the Arithmetic
without parameters\textup.
\end{TheoremB}

\begin{TheoremC}
If\/ $G$ is a group as in Theorem\/ \textup{A}\textup,
then the elementary theory of\/ $G$ is
hereditarily undecidable\textup.
\end{TheoremC}

In particular, all the groups of Thompson and Higman
interpret the Arithmetic with parameters,
while Thompson's group $F$ also interprets it without parameters,
and the elementary theories of all these groups
are hereditarily undecidable.

Theorems A and B are proved in
Section \ref{section:interpretationarithmetic}.
To the best of our knowledge, the interpretation constructed in the proof
of Theorem B is entirely original.
Theorem C is proved in Section \ref{section:undecidability}
as a corollary of Theorem A.
In the Appendix we show that every element of the derived subgroup of
$\F(r,\Lambda,A)$ is the product of two commutators,
and hence the derived subgroup is definable in $\F(r,\Lambda,A)$.

The main idea of the proof of Theorem A is,
as in \cite{BardakovTolstykh:2007:iaTgF},
to find in $G$ a definable subgroup isomorphic to
the restricted wreath product $\mathbb Z\wr\mathbb Z$,
because it is known that the latter group interprets the Arithmetic.
Note that, as opposed to $\mathbb Z\wr\mathbb Z$ and to the
groups that we study here, neither abelian groups,
nor virtually abelian,
nor free groups, nor torsion-free hyperbolic ones
can interpret the Arithmetic because their elementary theories
are all \emph{stable},
while the elementary theory of every structure that interprets
the Arithmetic is ``strongly'' unstable.
\emph{Stability\/} is a fundamental concept of Model Theory.
The textbooks
\cite{Pillay:1983:ist,Poizat:1985:ctm-fr,Poizat:2000:cmt-eng}
are all excellent introductions to the subject.
The best sources for learning about \emph{stable groups},
\ie\ the groups whose elementary theories are stable, are,
in our opinion,
\cite{Poizat:1987:gs-fr,Poizat:2001:sg-eng,Wagner:2000:sg}.
A proof of the stability of abelian groups can be found in
\cite[Theorem 3.1]{Prest:1988:mtm}.
Every non-elementary torsion-free hyperbolic group is stable
according to a recent result of Zlil Sela
\cite{Sela:pp2006:dgg8s}.

\begin{figure}
\includegraphics{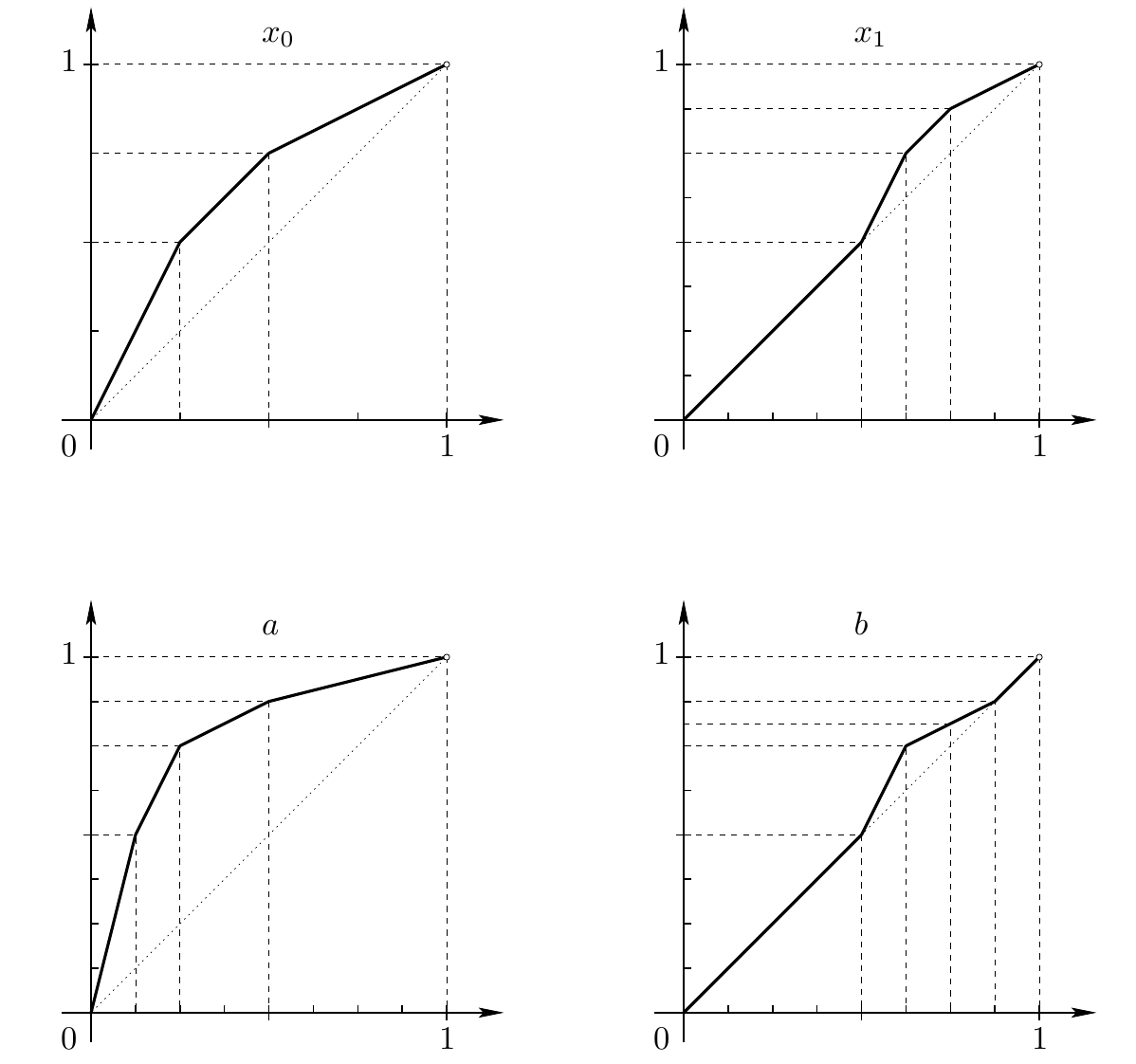}
\caption{Maps $x_0$, $x_1$, $a=x_0^2$,
and $b=x_1x_0^{-1}x_1^{-1}x_0$.}
\label{figure:1}
\end{figure}

A definable subgroup of $F$ isomorphic to $\mathbb Z\wr\mathbb Z$
was chosen by Bardakov and Tolstykh
\cite{BardakovTolstykh:2007:iaTgF} as follows.
Let $x_0$ and $x_1$ be the ``standard'' generators of $F$,
and let $a=x_0^2$, $b=x_1x_0^{-1}x_1^{-1}x_0$
(see Figure \ref{figure:1}).
One can verify without major difficulty that
$\langle a,b\rangle
=\langle b\rangle\wr\langle a\rangle\cong\mathbb Z\wr\mathbb Z$.
The centralizer of $x_0$ in $F$ coincides with the subgroup
generated by $x_0$.
Consequently, the subgroup $\langle a\rangle$ is definable
in $F$ by a formula with the parameter $x_0$.
Then it is shown that the centralizer of the subset
$\{\,a^{-k}ba^{k}\mid k\in\mathbb Z\,\}$ coincides with
the subgroup $\langle\,a^{-k}ba^{k}\mid k\in\mathbb Z\,\rangle$.
As the subset $\{\,a^{-k}ba^{k}\mid k\in\mathbb Z\,\}$
is clearly definable with the parameters $x_0$ and $x_1$,
so is the subgroup
$\langle\,a^{-k}ba^{k}\mid k\in\mathbb Z\,\rangle$.
Thus the subgroup $\langle a,b\rangle\cong\mathbb Z\wr\mathbb Z$
is definable in $F$ with parameters.
In the proof of Theorem A, we follow a similar approach for the group
$\F(r,\Lambda,A)$.


\section{Generalities}
\label{section:generalities}

In this section we present some basic definitions and facts.

\subsection{Permutations and piecewise affine maps}
\label{subsection:generalities.permutations}

\begin{definition}
A bijection of a set onto itself shall be called a \emph{permutation\/} of
this set.
A map $f$ shall be said to \emph{permute\/} a set $S$ if
the restriction $f|_S$ is a permutation of~$S$.
\end{definition}

\begin{definition}
Let $S$ be a set and $f$ be a bijection of $S$ onto itself.
We shall call the \emph{support\/} of $f$, denoted $\supp(f)$,
the complement in $S$ of the set of \emph{fixed points\/} of $f$,
denoted $\fix(f)$.
\end{definition}

Customarily, in the context of study of Thompson and Higman's groups,
\emph{all maps act on the right}.
We adopt the same convention in this article; for example:
$(\alpha)(xy)=((\alpha)x)y$ if
$x$ and $y$ are permutations of a set $S$, and
$\alpha\in S$.

We are going to write $X^f$, or sometimes $(X)f$, to denote
the image of the set $X$ under the map~$f$.

Lemmas \ref{lemma:01.(2.1)}, \ref{lemma:02.(2.2)},
and \ref{lemma:03.(2.3)} are obvious:

\begin{lemma}
\label{lemma:01.(2.1)}
Two permutations of the same set commute
if their supports are disjoint\textup.
\end{lemma}

\begin{lemma}
\label{lemma:02.(2.2)}
Let\/ $f$ and\/ $g$ be two permutations of the same set\textup.
Then
\begin{align*}
\fix(g^{-1}fg)&=(\fix(f))g,\\
\supp(g^{-1}fg)&=(\supp(f))g.
\end{align*}
\end{lemma}

\begin{lemma}
\label{lemma:03.(2.3)}
Let\/ $f$ and\/ $g$ be two commuting permutations of the same set\textup.
Then\/ $g$ permutes each of the sets\/
$\fix(f)$ and\/ $\supp(f)$\textup.
\end{lemma}

\begin{lemma}
\label{lemma:04.(2.4)}
Let\/ $I$ be an interval in\/ $\mathbb R$\textup,
$f$ an increasing permutation of\/ $I$\textup, and\/
$n\in\mathbb Z\setminus\{0\}$\textup.
Then\/ $\supp(f^n)=\supp(f)$\textup.
\end{lemma}
\begin{proof}
Clearly $\fix(f)\subset\fix(f^n)$ and $\supp(f)\supset\supp(f^n)$.
Consider an arbitrarily chosen $\alpha\in\supp(f)$.
Without loss of generality, suppose that $(\alpha)f>\alpha$.
Then
$$
\alpha<(\alpha)f<(\alpha)f^2<\dotsb<(\alpha)f^n,
$$
and hence $\alpha\in\supp(f^n)$.
\end{proof}

\begin{lemma}
\label{lemma:05.(2.5)}
Let\/ $I$ be a compact interval in\/ $\mathbb R$\textup,
$f$ be an increasing permutation of\/ $I$\textup, and\/
$\alpha\in I$\textup.
Then
$$
\lim_{n\to+\infty}(\alpha)f^n\in\fix(f).
$$
\end{lemma}
\begin{proof}
As $f$ preserves the order, the sequence
$((\alpha)f^n)_{n=0,1,\dotsc}$ is monotone, and hence the limit exists
and belongs to $I$.
If $\beta=\lim_{n\to+\infty}(\alpha)f^n$, then
$(\beta)f=\beta$ by continuity.
\end{proof}

\begin{definition}
Let $f$ be a map and $\alpha$ be a real number.
The map $f$ shall be called
\emph{affine to the right of\/ $\alpha$} if there exists
$\beta>\alpha$ such that the restriction $f|_{(\alpha,\beta)}$ is
an affine map $(\alpha,\beta)\to\mathbb R$.
The map $f$ shall be called
\emph{affine to the left of\/ $\alpha$} if there exists
$\beta<\alpha$ such that $f|_{(\beta,\alpha)}$ is
an affine map $(\beta,\alpha)\to\mathbb R$.
\end{definition}

\begin{definition}
For all $\alpha,\beta\in\mathbb R$ such that $\alpha<\beta$, and
for every map $f$ such that $f|_{(\alpha,\beta)}$ is an
affine map $(\alpha,\beta)\to\mathbb R$,
we shall denote the slope of $f|_{(\alpha,\beta)}$
by $(\alpha)f^{\prime+}$
and also by $(\beta)f^{\prime-}$.
The number $(\alpha)f^{\prime+}$ shall be called
\emph{the slope of\/ $f$ to the right of\/ $\alpha$}, and
$(\beta)f^{\prime-}$ shall be called
\emph{the slope of\/ $f$ to the left of\/~$\beta$}.
\end{definition}

The following lemma can be proved the same way in which one determines
the derivative of a composite function.

\begin{lemma}
\label{lemma:06.(2.6)}
Let\/ $f$ and\/ $g$ be two maps between subsets of\/
$\mathbb R$\textup, and let $\alpha\in\mathbb R$\textup.
\begin{enumerate}
\item
	If\/ $f$ is affine to the right of\/ $\alpha$ with\/
	$(\alpha)f^{\prime+}>0$\textup,
	$f$ is right-continuous at\/ $\alpha$\textup, and\/
	$g$ is affine to the right of\/ $(\alpha)f$\textup, then
	$$
	(\alpha)(fg)^{\prime+}=(\alpha)f^{\prime+}\cdot((\alpha)f)g^{\prime+}.
	$$
\item
	If\/ $f$ is affine to the left of\/ $\alpha$ with\/
	$(\alpha)f^{\prime-}>0$\textup,
	$f$ is left-continuous at\/ $\alpha$\textup, and\/
	$g$ is affine to the left of\/ $(\alpha)f$\textup, then
	$$
	(\alpha)(fg)^{\prime-}=(\alpha)f^{\prime-}\cdot((\alpha)f)g^{\prime-}.
	$$
\end{enumerate}
\end{lemma}

\subsection{Groups under consideration}
\label{subsection:generalities.groups}

As in the Introduction,
choose a positive real number $r$, a subgroup
$\Lambda$ of the multiplicative group $\mathbb R^*_+$, and
a submodule $A$ of the $\Lambda$-module $\mathbb R$
such that $r\in A$.
Call such a triple $(r,\Lambda,A)$ \emph{admissible}.
For an admissible triple $(r,\Lambda,A)$,
define the groups
$\F(r,\Lambda,A)$, $\T(r,\Lambda,A)$, and $\V(r,\Lambda,A)$
as in the Introduction.
We are going to treat Thompson's groups $F$, $T$, and $V$
as particular cases, so we pose
$$
F=\F(1,\langle 2\rangle,\mathbb Z[{\textstyle\frac{1}{2}}]),\qquad
T=\T(1,\langle 2\rangle,\mathbb Z[{\textstyle\frac{1}{2}}]),\qquad
V=\V(1,\langle 2\rangle,\mathbb Z[{\textstyle\frac{1}{2}}]).
$$

\begin{remark}
\label{remark:01.(2.7)}
Every element of $\F(r,\mathbb R^*_+,\mathbb R)$ can be extended in a
unique way to a homeomorphism $[0,r]\to[0,r]$
with respect to the usual topology.
\end{remark}

\begin{remark}
\label{remark:02.(2.8)}
The conjugation of elements of $\V(1,\mathbb R^*_+,\mathbb R)$
by the linear map of multiplication by $r$ is
an isomorphism between $\V(1,\mathbb R^*_+,\mathbb R)$
and $\V(r,\mathbb R^*_+,\mathbb R)$, which maps
$\V(1,\Lambda,Ar^{-1})$ onto $\V(r,\Lambda,A)$,
$\T(1,\Lambda,Ar^{-1})$ onto $\T(r,\Lambda,A)$, and
$\F(1,\Lambda,Ar^{-1})$ onto $\F(r,\Lambda,A)$,
where $Ar^{-1}=\{\,\alpha/r\mid \alpha\in A\,\}$.
\end{remark}

It is convenient to distinguish in $\F(r,\Lambda,A)$
the following subsemigroups:
\begin{itemize}
\item
	denote by $\F^\uparrow(r,\Lambda,A)$ the semigroup of all the elements
	$x\in\F(r,\Lambda,A)$ such that
	$(\alpha)x\ge\alpha$ for all $\alpha\in[0,r)$, and
\item
	denote by $\F^\downarrow(r,\Lambda,A)$
	the semigroup of all the elements
	$x\in\F(r,\Lambda,A)$ such that
	$(\alpha)x\le\alpha$ for all $\alpha\in[0,r)$.
\end{itemize}

We shall introduce other useful notations:
if $S\subset[0,r)$, then
\begin{align*}
\F_S(r,\Lambda,A)
&=\{\,x\in\F(r,\Lambda,A)\mid\supp(x)\subset S\,\},\\
\T_S(r,\Lambda,A)
&=\{\,x\in\T(r,\Lambda,A)\mid\supp(x)\subset S\,\},\\
\V_S(r,\Lambda,A)
&=\{\,x\in\V(r,\Lambda,A)\mid\supp(x)\subset S\,\}.
\end{align*}
The semigroups
$\F^\uparrow_S(r,\Lambda,A)$ and $\F^\downarrow_S(r,\Lambda,A)$
are defined similarly.

For the rest of this article, we fix $(r,\Lambda,A)$,
and moreover, we assume that
\emph{$\Lambda$ is nontrivial\/}: $\Lambda\ne\{1\}$.
To simplify the notation, we abbreviate:
$$
\F=\F(r,\Lambda,A),\qquad
\T=\T(r,\Lambda,A),\qquad
\V=\V(r,\Lambda,A).
$$
We are also going to write $\F^\uparrow$ instead of
$\F^\uparrow(r,\Lambda,A)$, etc.

\subsection{Theories and models}
\label{subsection:generalities.modeltheory}

In this article we talk about \emph{structures\/} in the sense of
Model Theory
(or in the sense of Universal Algebra,
up to a few linguistic distinctions).
When the terms ``formula,'' ``sentence,'' and ``theory'' are used
in the formal sense, they always denote first-order formulae,
sentences, and theories.
The terms ``theory,'' ``elementary theory,'' and ``first-order theory''
are to be used as synonyms.
``A model'' and ``a structure'' will usually be synonyms as well.
Except when indicated otherwise, the formulae are
\emph{without parameters}.

A structure $M$ of \emph{signature\/} $\Sigma$,
also called \emph{$\Sigma$-structure\/},
is a \emph{model\/}
of a set of sentences $S$, this fact being denoted $M\models S$, if
$M$ \emph{satisfies\/} every sentence $\alpha$ of $S$,
which be denoted $M\models\alpha$.
A sentence $\alpha$ is called a \emph{consequence\/}
of a set of sentences $S$ in a signature $\Sigma$,
which be denoted $S\vdash_\Sigma\alpha$,
or $S\vdash\alpha$ in the case when $\Sigma$ is well understood, if
every $\Sigma$-model of $S$ is also a model of $\alpha$.
A set of sentences is \emph{consistent\/}
in a signature $\Sigma$ if
it has a $\Sigma$-model.%
\footnote{The signature $\Sigma$ in this definition
is of minor importance except in the case when
$S$ has an empty model,
which would imply in particular that no element of 
$S$ contain any constant symbols.}
A set of sentences is \emph{deductively closed\/}
in a signature $\Sigma$ if
it contains all its consequences in $\Sigma$.
A \emph{theory of signature\/ $\Sigma$}, also called
\emph{$\Sigma$-theory}, is a consistent and deductively closed
in $\Sigma$ set of sentences.
If $T$ is a $\Sigma$-theory, then $T\vdash_\Sigma\alpha$ is equivalent to
$\alpha\in T$.
The \emph{theory of a structure\/ $M$}, denoted $\Th(M)$,
is the set of all the sentences in the signature of $M$
satisfied by $M$.
A theory is called \emph{complete\/} if
it is the theory of a structure.
The class of all the $\Sigma$-models of a set of
$\Sigma$-sentences $S$ is denoted $\Mod_\Sigma(S)$.

We are going to use implicitly \emph{the Compactness Theorem},
which guaranties that a consequence of a set of sentences is always
a consequence, in the same signature, of its finite subset.
We recommend any of the textbooks
\cite{Hodges:1993:mt,Hodges:1997:smt,
Poizat:1985:ctm-fr,Poizat:2000:cmt-eng,Rothmaler:2000:imt-eng}
for references on general results of Model Theory.

If $M$ is a structure,
a \emph{set definable in\/ $M$} in general is
a subset of $M^n$, where $n\in\mathbb N$,
definable by a first-order formula in the signature of $M$
and possibly with \emph{parameters\/} from $M$.
The \emph{parameters\/} are new constant symbols added
to the language and interpreted by elements of $M$.
(Usually, to name an element $a\in M$,
one uses $a$ itself as a parameter.)
For example, in a group, the centralizer of every element $g$ is
definable by the formula $\phi(x)=\ulcorner gx=xg\urcorner$
with $g$ as a parameter,
but in general there is no reason for that centralizer
to be definable by a formula without parameters in the pure-group
signature, which only contains a single binary function symbol 
$\ulcorner\cdot\urcorner$ to denote the group operation
$(x,y)\mapsto x\cdot y$.
To be explicit, a set is said to be
``definable \emph{with\/}'' or ``\emph{without\/} parameters''
according to whether or not parameters are allowed in its definition.
However, we are not going to specify in which structure a given set
is definable, except if there are several equally natural choices
in the context.
Of course it is possible as well to talk about definability
of relations and operations.

If $f$ is a map $A\to B$, and $n$ is a positive integer,
we denote by $f^{\underline{n}}$ the map $A^n\to B^n$
induced by $f$.
We are going to slightly abuse the notation by assuming that if
the domain of $f$ is a subset of $A^m$, then the domain
of $f^{\underline{n}}$ is naturally identified with a subset of $A^{mn}$.
We shall call the \emph{$f$-preimage\/} of a given set
its preimage under $f^{\underline{n}}$
when the choice of $n$ is clear
(so not necessarily under $f$ itself).

Consider two structures,
$M$ of signature $\Sigma$ and $N$ of signature $\Gamma$.

\begin{definition}
We call an \emph{interpretation of\/ $M$ in\/ $N$ with parameters\/}
a pair $(n,f)$ where
$n\in\mathbb N$ and $f$ is a surjective map of a subset of
$N^n$ onto $M$
such that for every set $X$ definable in $M$ without parameters,
the $f$-preimage of $X$ is definable (in $N$)
with (possibly) parameters.
An interpretation $(n,f)$ with parameters is called an
\emph{interpretation without parameters\/} if
the $f$-preimage of every set definable without parameters is
definable without parameters as well.
\end{definition}

See \cite[Chapter 5]{Hodges:1993:mt} for detailed explanation
of interpretability and related notions.

In what follows, in accordance with conventions of Model Theory,
the terms ``definable,'' ``$0$-definable,''
``intepretation,'' ``$0$-interpretation'' shall be used to denote,
respectively,
``definable with parameters,''
``definable without parameters,''
``interpretation with parameters,'' and
``interpretation without parameters.''
Furthermore, as in our case the value of $n$ for an interpretation 
$(n,f)$ under consideration will be often either well understood,
or hardly important, we shall simplify the notation and 
call $f$ itself an interpretation.
To indicate that $(n,f)$ is an interpretation of $M$ in $N$,
we shall write
either $(n,f)\colon M\rightsquigarrow N$,
or $f\colon M\rightsquigarrow N$.%
\footnote{The arrow here is pointing in the opposite direction
than that in the notation of \cite{AhlbrandtZiegler:1986:qfatct}.}
To indicate that $f$ is a $0$-interpretation of $M$ in $N$,
we shall write
$f\colon M\stackrel{\scriptscriptstyle\varnothing}{\rightsquigarrow}N$.
To indicate that $M$ is interpretable or $0$-interpretable
in $N$, we write $M\rightsquigarrow N$ or
$M\stackrel{\scriptscriptstyle\varnothing}{\rightsquigarrow}N$,
respectively.

\begin{remark}
\label{remark:03.(2.9)}
If $f\colon M\rightsquigarrow N$, then
the $f$-preimage of every set definable in $M$ is also definable (in $N$).
\end{remark}

\begin{remark}
\label{remark:04.(2.10)}
Let $n\in\mathbb N$ and $B\subset N^n$.
Then a surjective map $f$ of $B$ onto $M$ is an interpretation
of $M$ in $N$ if and only if
\begin{enumerate}
\item
	the domain $B$ is definable,
\item
	the equivalence relation on $B$ induced by $f$
	(the \emph{kernel\/} of $f$) is definable, and
\item
	for every relation, operation, and constant of the structure $M$
	(named by a symbol of $\Sigma$),
	the $f$-preimage of its graph is definable.
\end{enumerate}
The map $f$ is a $0$-interpretation if and only if
all these sets are $0$-definable.
\end{remark}

\begin{remark}
\label{remark:05.(2.11)}
If $L$, $M$, and $N$ are three structures, and
$(m,f)\colon L\rightsquigarrow M$ and $(n,g)\colon M\rightsquigarrow N$,
then $(mn,g^{\underline{m}}f)\colon L\rightsquigarrow N$.
If moreover $f$ and $g$ are $0$-interpretations, then
so is $g^{\underline{m}}f$.
\end{remark}

\begin{definition}
[see \cite{AhlbrandtZiegler:1986:qfatct} et
{\cite[Section 5.4(c)]{Hodges:1993:mt}}]
Two structures $M$ and $N$ are said to be \emph{bi-interpretable\/} if
there exist two interpretations
$(m,f)\colon M\rightsquigarrow N$ and $(n,g)\colon N\rightsquigarrow M$
such that
the map $g^{\underline{m}}f$ is definable in $M$,
and $f^{\underline{n}}g$ is definable in $N$.
The interpretations $(m,f)$ and $(n,g)$ in this case are called
\emph{bi-interpretations}.
\end{definition}

\subsection{Decidability}
\label{subsection:generalities.decidability}

Let $A$ be a finite set viewed as an alphabet,
and denote $A^*$ the set of all finite words in $A$.
We say that a set $X\subset A^*$ is \emph{recursive\/}
or \emph{decidable\/} if
there exists an \emph{algorithm\/} which for every input $w\in A^*$
answers the question whether $w\in X$.
A map $f\colon X\to A^*$ is said to be \emph{computable\/} if
there exists an algorithm which computes $f(w)$ for every input $w\in X$,
and which never stops for any input $w\notin X$.

Usually a set is said to be recursive or non-recursive,
while a theory is said to be decidable or undecidable.

In the rest of this section, let $\Sigma$ be an arbitrary
finite signature.

\begin{definition}
A $\Sigma$-theory $T$
is said to be \emph{essentially undecidable\/} if
every $\Sigma$-theory containing $T$ is undecidable.
\end{definition}

\begin{remark}
\label{remark:06.(2.12)}
If $T$ is a theory of finite signature,
and a subset of $T$ is an essentially undecidable theory
(possibly of smaller signature),
then $T$ is undecidable, and even essentially undecidable.
\end{remark}

\begin{definition}
A $\Sigma$-theory $T$
is said to be \emph{hereditarily undecidable\/} if
every $\Sigma$-subtheory of $T$ is undecidable.
\end{definition}

\begin{lemma}[{\cite[Theorem 6]{Tarski:1971:gmpu}}]
\label{lemma:07.(2.13)}
If\/ $T$ is a theory of finite signature\textup,
and\/ $T$ has a finitely axiomatized essentially undecidable
subtheory\/
\textup(possibly of smaller signature\/\textup)\textup,
then\/ $T$ is hereditarily undecidable\textup.
\end{lemma}
\begin{proof}
Denote by $\Sigma$ the signature of $T$.
Let $S$ be a finitely axiomatized essentially undecidable
subtheory of $T$.
Choose a sentence $\theta\in S$ which axiomatizes~$S$.

Arguing by contradiction,
let $U$ be a decidable $\Sigma$-subtheory of $T$.
Let $R$ be the $\Sigma$-theory axiomatized (generated) by $U\cup S$.
Then
$$
R=\{\,\alpha\mid U,\theta\vdash_\Sigma\alpha\,\}
=\{\,\alpha\mid\ulcorner\theta\rightarrow\alpha\urcorner\in U\,\},
$$
and hence $R$ is decidable since so is $U$.
This contradicts the essential undecidability of $S$
(see Remark~\ref{remark:06.(2.12)}).
\end{proof}


\section{Two lemmas}
\label{section:twolemmas}

In this section we will prove two technical lemmas about
piecewise affine homeomorphisms of an interval.
These lemmas will be essential for the proof of Theorem A, namely
for showing that certain centralizers are preserved when passing
from $\F$ to $\V$, and thus being able to pass from an interpretation
of the Arithmetic in $\F$ to its interpretations in $\V$ and~$\T$.

\begin{lemma}
\label{lemma:08.(3.1)}
Let\/ $r\in\mathbb R^*_+$\textup.
Let\/ $z$ be a homeomorphism\/ $[0,r)\to[0,r)$ such that\/
$(\alpha)z>\alpha$ for all\/ $\alpha\in(0,r)$\textup.
Let\/ $f$ be a permutation of\/ $[0,r)$
such that\/\textup:
\begin{enumerate}
\item
	$f$ commutes with\/ $z$\textup,
\item
	$f$ is\/ \textup(right-\/\textup)continuous at\/ $0$\textup,
\item
	$f$ has a finite number of points of discontinuity\textup.
\end{enumerate}
Then\/ $f$ is continuous\/
\textup(with respect to the usual topology\/\textup)\textup.
\end{lemma}
\begin{proof}
Suppose that $f$ were not continuous.
Then let $\alpha$ be the least element of $[0,r)$
at which $f$ is not continuous.
Since $f$ is continuous at $0$, $\alpha\in(0,r)$.
Hence $(\alpha)z^{-1}<\alpha$ and $f$ is continuous at $(\alpha)z^{-1}$.
It follows that $f$ is continuous at $\alpha$ because
$f=z^{-1}fz$, where $z^{-1}$ and $z$ are continuous everywhere,
and $f$ is continuous at $(\alpha)z^{-1}$.
This gives a contradiction.
\end{proof}

\begin{lemma}
\label{lemma:09.(3.2)}
Let\/ $r\in\mathbb R^*_+$\textup.
Let\/ $Z$ be a set of homeomorphisms\/ $[0,r)\to[0,r)$
such that\/\textup:
\begin{enumerate}
\item
	$(\alpha)z\ge\alpha$\/ for all\/ $z\in Z$ and all\/
	$\alpha\in[0,r)$\textup,
\item
	$\supp(z)$ is an interval for every\/ $z\in Z$\textup,
\item
	$\bigcup_{z\in Z}\supp(z)$ is dense in\/ $(0,r)$\textup.
\end{enumerate}
Let\/ $f$ be a permutation of\/ $[0,r)$ such that\/\textup:
\begin{enumerate}
\item
	$f$ commutes with every\/ $z\in Z$\textup,
\item
	$f$ is right-continuous at every point of\/ $[0,r)$\textup,
\item
	$f$ has a finite number of points of discontinuity\textup.
\end{enumerate}
Then\/ $f$ is continuous\textup.
\end{lemma}
\begin{proof}
We shall show first that $f$ is increasing on each of its intervals
of continuity.

We shall say that an interval $I$ is \emph{right-closed\/} if
$\sup I\in I$,
and that it is \emph{left-closed\/} if
$\inf I\in I$.
If an interval is not right- or left-closed, we shall say that is
\emph{right-\/} or \emph{left-open}, respectively.

As the number of points of discontinuity of $f$ is finite,
the number of its maximal intervals of continuity is finite as well.
Since $f$ is right-continuous everywhere in its domain $[0,r)$,
every maximal interval of continuity is left-closed and right-open.
Therefore the image under $f$ of every maximal interval of continuity $I$
is left-closed if $f$ increases on $I$ and right-closed
if $f$ decreases on~$I$.

Clearly $[0,r)$ is the disjoint union of the images of the maximal
intervals of continuity of $f$.
Since every such image is either a right- or a left-open interval,
the only possibility is that they all are left-closed and right-open.
Hence $f$ is increasing on each of its intervals of continuity.

Suppose now that $f$ were not continuous.

By Lemma \ref{lemma:03.(2.3)}, for every $z\in Z$,
$f$ permutes the set $\supp(z)$.
We shall show that $f$ is continuous on each of these intervals.
Arguing by contradiction, take $z\in Z$ such that
$f$ is not continuous on $\supp(z)$.
Let $\gamma$ be the least element of $\supp(z)$ at which
$f$ is not continuous.
Then $f$ is continuous at $\gamma$ because
$f=z^{-1}fz$, where $z$ is a homeomorphism,
and $f$ is continuous at $(\gamma)z^{-1}$, since $(\gamma)z^{-1}<\gamma$.
This gives a contradiction.

For every $z\in Z$, $f$ is increasing on $\supp(z)$
since it is continuous there.
Therefore for every $z\in Z$,
$$
\lim_{\alpha\to(\sup\supp(z))^{-}}(\alpha)f=\sup\supp(z).
$$
By right continuity, $(\inf\supp(z))f=\inf\supp(z)$
for all $z\in Z$.

Let $S=\bigcup_{z\in Z}\supp(z)$ and
$L=\{\,\inf\supp(z)\mid z\in Z\,\}$.
Then $S$ is open and dense in $[0,r)$, $f|_S$ is continuous,
and $f|_L=\id_L$.

Consider an arbitrary $\alpha\in[0,r)\setminus S$.
It is easy to see that for any $\beta\in(\alpha,r)$,
$L\cap[\alpha,\beta)$ is not empty.
Hence $(\alpha)f=\alpha$ by right continuity at $\alpha$.
We have shown that
$f|_{[0,r)\setminus S}=\id_{[0,r)\setminus S}$.

The map $f$ is increasing.
Indeed, if $f$ is a map from a linearly ordered set into itself,
and this set is covered by intervals such that $f$ sends each of them
into itself in a strictly increasing way, then $f$ is strictly increasing.
In our case we have:
$$
[0,r)=\bigcup_{z\in Z}\supp(z)
\cup\bigcup_{\alpha\in[0,r)\setminus S}[\alpha,\alpha].
$$

Therefore $f$ is continuous as an increasing surjection
of a subset of $\mathbb R$ onto an interval of~$\mathbb R$.
\end{proof}

These last two lemmas already suffice to show, using the result
of Bardakov and Tolstykh, that
the definable subgroup of $F$ isomorphic to $\mathbb Z\wr\mathbb Z$
used in \cite{BardakovTolstykh:2007:iaTgF} is definable in
$T$ and $V$ as well.

\begin{proposition}
\label{proposition:01.(3.3)}
Let\/ $a$ and\/ $b$ be the elements of\/ $F$ shown on
Figure\/ \textup{\ref{figure:1}.}
Then
$$
\langle a,b\rangle
=\langle b\rangle\wr\langle a\rangle\cong\mathbb Z\wr\mathbb Z,
$$
and the subgroup\/ $\langle a,b\rangle$ is definable
in\/ $F$\textup, in\/ $T$\textup, and in\/ $V$
by the same first-order formula with parameters\textup.
\end{proposition}
\begin{proof}
We consider the same elements
$x_0$, $x_1$, $a=x_0^2$, and $b=x_1x_0^{-1}x_1^{-1}x_0$ of $F$
as in \cite{BardakovTolstykh:2007:iaTgF},
see Figure \ref{figure:1}.
It is shown in \cite{BardakovTolstykh:2007:iaTgF} that:
\begin{enumerate}
\item
	$\langle a,b\rangle
	=\langle b\rangle\wr\langle a\rangle\cong\mathbb Z\wr\mathbb Z$,
\item
	$C_F(x_0)=\langle x_0\rangle$,
\item
	$C_F(\{\,a^{-k}ba^{k}\mid k\in\mathbb Z\,\})
	=\langle\,a^{-k}ba^{k}\mid k\in\mathbb Z\,\rangle$.
\end{enumerate}
In particular, $\langle a,b\rangle$ is the semi-direct product of
$C_F(\{\,a^{-k}ba^{k}\mid k\in\mathbb Z\,\})$ and $\langle a\rangle$,
where $\langle a\rangle=\{\,x^2\mid x\in C_F(x_0)\,\}$.

Define $\alpha_k=2^{-1+2k}$ for $k=-1,-2,-3,\dotsc$, and
$\alpha_k=1-2^{-1-2k}$ for $k=0,1,2,\dotsc$.
Then
$$
0<\dotsb<\alpha_{-2}<\alpha_{-1}<\alpha_0<\alpha_1<\alpha_2<\dotsb<1.
$$

A direct calculation facilitated by Lemma \ref{lemma:02.(2.2)}
shows that:
\begin{enumerate}
\item
	$\supp(x_0)=(0,1)$;
\item
	$(\alpha)x_0\ge\alpha$ for all $\alpha\in[0,1)$;
\item
	$\supp(a^{-k}ba^k)=(\alpha_k,\alpha_{k+1})$ for all $k\in\mathbb Z$;
\item
	$(\alpha)a^{-k}ba^k\ge\alpha$
	for all $\alpha\in[0,1)$ and all $k\in\mathbb Z$.
\end{enumerate}

Choosing $x_0$ as $z$ in Lemma \ref{lemma:08.(3.1)},
we conclude that every element of $C_V(x_0)$ is continuous.
Setting $Z=\{\,a^{-k}ba^{k}\mid k\in\mathbb Z\,\}$
in Lemma \ref{lemma:09.(3.2)},
we conclude that every element of
$C_V(\{\,a^{-k}ba^{k}\mid k\in\mathbb Z\,\})$ is continuous.
Since an element of $V$ belongs to $F$ if and only if it is continuous,
the centralizer of the element $x_0$ and the centralizer of the set
$\{\,a^{-k}ba^{k}\mid k\in\mathbb Z\,\}$
are preserved when passing from $F$ to $V$.
Hence $\langle a,b\rangle$ is definable in $F$, $T$, and $V$
by the same first-order formula with parameters.
\end{proof}


\section{Definable copies of $\mathbb Z\wr\mathbb Z$}
\label{section:copiesZwrZ}

In this section we show that all the groups
$\F$, $\T$, and $\V$ have definable subgroups isomorphic to
$\mathbb Z\wr\mathbb Z$.
\begin{lemma}
\label{lemma:10.(4.1)}
Let\/ $\alpha,\beta\in[0,r]\cap A$ be such that\/ $\alpha<\beta$\textup,
and denote\/ $I=(\alpha,\beta)$\textup.
Let\/ $x\in\F_I$ be such that\/
$\fix(x)\cap I\cap A=\varnothing$\textup.
Let\/ $\phi\colon\F_I\to\Lambda$
be the map\/ $y\mapsto(\alpha)y^{\prime+}$\textup.
Let\/ $C$ be the centralizer of\/ $x$ in\/
$\F_I$\textup.
Then\/ $\phi$ is a homomorphism\textup, and
its restriction to\/ $C$ is injective\textup.
The same holds for\/
$\phi\colon \F_I\to\Lambda,\ y\mapsto(\beta)y^{\prime-}$\textup.
\end{lemma}
\begin{proof}
We shall only consider the case of
$\phi\colon\F_I\to\Lambda,\ y\mapsto(\alpha)y^{\prime+}$, because
the case of $y\mapsto(\beta)y^{\prime-}$ is analogous.
It is easy to verify that $\phi$ is a homomorphism
(see Lemma \ref{lemma:06.(2.6)}).
Suppose that $\phi$ is not injective on~$C$.

Let $y\in C$ be such that $\phi(y)=1$ but $y\ne\id$.
Let $\gamma\in(\alpha,\beta)$ be such that
$$
y|_{[0,\gamma]}=\id_{[0,\gamma]}\quad\text{but}\quad
(\gamma)y^{\prime+}\ne1.
$$
Then $\gamma\in A$, and hence $(\gamma)x\ne\gamma$.
Without loss of generality, assume that $(\gamma)x>\gamma$,
because otherwise $(\gamma)x^{-1}>\gamma$, and one can use
$x^{-1}$ instead of $x$.
Then
$$
[0,(\gamma)x]=[0,\gamma]^x\subset\fix(y),
$$
see Lemma \ref{lemma:03.(2.3)},
and hence $(\gamma)y^{\prime+}=1$.
This gives a contradiction.
\end{proof}

\begin{lemma}
\label{lemma:11.(4.2)}
The centralizer\/ $C$ in
Lemma\/ \textup{\ref{lemma:10.(4.1)}} is cyclic\textup.
\end{lemma}

This lemma follows from
the description of centralizers in $\F(r,\mathbb R^*_+,\mathbb R)$
obtained by Matthew Brin and Craig Squier
\cite{BrinSquier:2001:pcrcgplhrl},
but for the reader's convenience we prefer to provide a self-contained
proof.

\begin{proof}[Proof of Lemma\/ \textup{\ref{lemma:11.(4.2)}}]
First of all, note that if $\Lambda$ is cyclic itself,
the conclusion of this lemma is an obvious corollary of
Lemma~\ref{lemma:10.(4.1)}.

Let $\alpha$, $\beta$, $I$, $x$, $\phi$, and $C$
be such as in Lemma \ref{lemma:10.(4.1)}.
Without loss of generality, assume that $(\alpha)x^{\prime+}>1$.

We are going to use the fact that a multiplicative subgroup of
$\mathbb R^*_+$ is either cyclic
(the trivial subgroup is considered cyclic),
or dense in $\mathbb R^*_+$ with respect to the usual topology.

Let $\Gamma$ be the image of the group $C$ under the homomorphism
$\phi\colon \F_I\to\Lambda$.
By Lemma \ref{lemma:10.(4.1)}, $\phi$ is injective,
and hence $C\cong\Gamma$.
It remains to show that $\Gamma$ is not dense in $\mathbb R^*_+$.

Observe that $\fix(x)\cap I$ is a finite set, and that
$$
\fix(y)\cap I=\fix(x)\cap I\quad\text{for all}\quad
y\in C\setminus\{\id\}.
$$
Indeed, it is clear that $\fix(x)\cap I$ is finite because
$\fix(x)\cap I\cap A$ is empty and $A$ is dense in $\mathbb R$.
Consider now an arbitrary $y\in C\setminus\{\id\}$.
By the injectivity of $\phi$, $(\alpha)y^{\prime+}\ne1$.
If $\fix(y)\cap I\cap A$ were nonempty, then it would have
a least element $\gamma$,
and this $\gamma$ would be fixed by $x$: $(\gamma)x=x$.
This would contradict $\fix(x)\cap I\cap A=\varnothing$,
hence $\fix(y)\cap I\cap A=\varnothing$
and $\fix(y)\cap I$ is finite.
As $x$ and $y$ commute and each of them permutes $I$,
it follows that $x$ permutes $\fix(y)\cap I$
and that $y$ permutes $\fix(x)\cap I$.
Since these sets are finite, and because
$x$ and $y$ preserve the order, it follows that
$\fix(y)\cap I\subset\fix(x)$,
$\fix(x)\cap I\subset\fix(y)$,
an hence $\fix(y)\cap I=\fix(x)\cap I$.

Denote
$$
\beta_0=\min((\fix(x)\cap I)\cup\{\beta\});
$$
it exists but does not belong to $A$ unless $\beta_0=\beta$.
Then $(\gamma)x>\gamma$ for all $\gamma\in(\alpha,\beta_0)$,
since $(\alpha)x^{\prime+}>1$.

Choose $\alpha_1,\beta_1\in(\alpha,\beta_0)$ such that
$x^{-1}$ be affine on
$[\alpha,\alpha_1]$ and $x$ be affine on $[\beta_1,\beta_0)$.
Then $x$ is also affine on
$[\alpha,(\alpha_1)x^{-1}]$,
and $x^{-1}$ is affine on $[(\beta_1)x,\beta_0)$.

We shall show now that if
$y\in C$ and $(\alpha)y^{\prime+}>1$, then $y$ is affine on
$[\alpha,(\alpha_1)y^{-1}]$ and on $[\beta_1,\beta_0)$.
Consider one such $y$.
Thus $(\gamma)y>\gamma$ for all $\gamma\in(\alpha,\beta_0)$.
Then for every $\gamma\in(\alpha,\alpha_1]$,
$$
x^{-1}|_{[\alpha,\alpha_1]}\cdot y^{-1}|_{[\alpha,(\gamma)x^{-1}]}\cdot
x|_{[\alpha,(\alpha_1)x^{-1}]}
=y^{-1}|_{[\alpha,\gamma]},
$$
and for every $\gamma\in[\beta_1,\beta_0)$,
$$
x|_{[\beta_1,\beta_0)}\cdot y|_{[(\gamma)x,\beta_0)}\cdot
x^{-1}|_{[(\beta_1)x,\beta_0)}
=y|_{[\gamma,\beta_0)}.
$$
These obvious equalities imply that:
\begin{enumerate}
\item
	for every $\gamma\in(\alpha,\alpha_1]$,
	$y^{-1}$ is affine on $[\alpha,\gamma]$ if it is
	affine on $[\alpha,(\gamma)x^{-1}]$,
\item
	for every $\gamma\in[\beta_1,\beta_0)$,
	$y$ is affine on $[\gamma,\beta_0)$ if it is
	affine on $[(\gamma)x,\beta_0)$.
\end{enumerate}
This is possible only if $y^{-1}$ is affine on $[\alpha,\alpha_1]$,
and $y$ is affine on $[\beta_1,\beta_0)$.

Clearly $\lim_{n\to+\infty}(\gamma)x^{n}=\beta_0$
for all $\gamma\in(\alpha,\beta_0)$
(see the proof of Lemma \ref{lemma:05.(2.5)}).
Choose an integer $n$ such that
$$
(\alpha_1)x^{n}>\beta_1,
\quad\text{and hence}\quad
(\beta_1)x^{-n}<\alpha_1.
$$
Denote $p=(\beta_1)(x^{-n})^{\prime+}$.
Then $p^{-1}=((\beta_1)x^{-n})(x^{n})^{\prime+}$.
Choose $\gamma$ in $((\beta_1)x^{-n},\alpha_1]$
such that $x^{n}$ be affine on
$[(\beta_1)x^{-n},\gamma]$ (with the slope $p^{-1}$).

Suppose that $y$ is an element of $C$ such that
$$
1<(y)\phi<\frac{\gamma-\alpha}{(\beta_1)x^{-n}-\alpha}.
$$
Then $y^{-1}$ is affine on
$[\alpha,\alpha_1]$, $y$ is affine on $[\beta_1,\beta_0)$,
and $(\beta_0)y^{\prime-}<1$
(because $(\gamma)y>\gamma$ for all $\gamma\in(\alpha,\beta_0)$).
Denote $q=(\alpha)y^{\prime+}=(y)\phi$.
Then
$$
(\beta_1)x^{-n}<\alpha+(\gamma-\alpha)q^{-1}=(\gamma)y^{-1}
\le\alpha+(\alpha_1-\alpha)q^{-1}=(\alpha_1)y^{-1},
$$
and therefore $(\beta_1)x^{-n}<(\beta_1)x^{-n}y<\gamma$.
Since $y=x^{-n}yx^{n}$,
\begin{align*}
(\beta_1)y^{\prime+}
&=(\beta_1)(x^{-n})^{\prime+}
\cdot((\beta_1)x^{-n})y^{\prime+}
\cdot((\beta_1)x^{-n}y)(x^{n})^{\prime+}\\
&=pqp^{-1}=q>1.
\end{align*}
As $(\beta_1)y^{\prime+}=(\beta_0)y^{\prime-}<1$,
this gives a contradiction, which means that
$$
\Gamma\cap\Bigl(1,\frac{\gamma-\alpha}{(\beta_1)x^{-n}-\alpha}\Bigr)
=\varnothing.
$$
\end{proof}

\begin{lemma}
\label{lemma:12.(4.3)}
Let\/ $\alpha,\beta\in[0,r]\cap A$ be such that\/ $\alpha<\beta$\textup.
Let\/ $p,q\in\Lambda$ be such that\/ $p>1>q$\textup.
Then there exists\/ $x\in\F^\uparrow$ such that\/
$\supp(x)=(\alpha,\beta)$\textup,
$(\alpha)x^{\prime+}=p$\textup, and\/ $(\beta)x^{\prime-}=q$\textup.
\end{lemma}
\begin{proof}
Let $s$ be an element of $\Lambda$ such that $(2+p+q)s\le1$.
Denote $l=\beta-\alpha$.
Consider two subdivisions of the interval $[\alpha,\beta]$:
the first one---into subintervals of lengths
$sl$, $qsl$, $(1-(2+p+q)s)l$, $psl$, and $sl$,
in this order,
and the second---into subintervals of lengths
$psl$, $sl$, $(1-(2+p+q)s)l$, $sl$, and $qsl$, in this order.
Let $x$ be the continuous map $[0,r]\to[0,r]$ which is
the identity on $[0,\alpha]\cup[\beta,r]$, and
which sends every interval of the first subdivision of $[\alpha,\beta]$
in the affine manner onto the corresponding interval of the second.
It is easy to verify that 
$\supp(x)=(\alpha,\beta)$,
$(\alpha)x^{\prime+}=p$, $(\beta)x^{\prime-}=q$,
and $x|_{[0,r)}\in\F^\uparrow$.
\end{proof}

\begin{figure}
\includegraphics{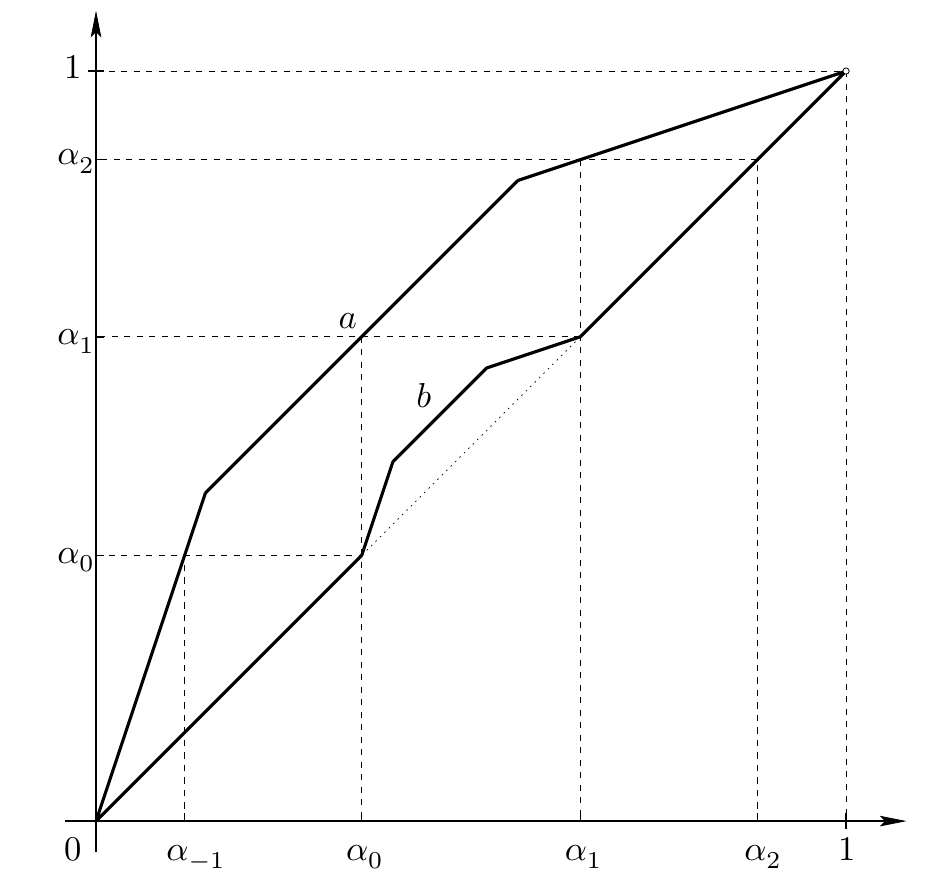}
\caption{An example of the maps $a$ and $b$.}
\label{figure:2}
\end{figure}

Choose $a\in\F^\uparrow$ such that $\supp(a)=(0,r)$.
Choose $\alpha_0\in(0,r)\cap A$ arbitrarily.
For every $k\in\mathbb Z$, define $\alpha_k=(\alpha_0)a^k$.
Observe that
\begin{gather*}
0<\dotsb<\alpha_{-2}<\alpha_{-1}<\alpha_0<\alpha_1<\alpha_2<\dotsb<r,\\
\text{and}\qquad
\lim_{n\to-\infty}\alpha_n=0,\quad\lim_{n\to+\infty}\alpha_n=r
\end{gather*}
(see Lemma \ref{lemma:05.(2.5)} and
Remark \ref{remark:01.(2.7)}).
Choose $b\in\F^\uparrow$ such that
$\supp(b)=(\alpha_0,\alpha_1)$
(see Figure \ref{figure:2}).
Such $a$ and $b$ exist by Lemma~\ref{lemma:12.(4.3)}.

\begin{lemma}
\label{lemma:13.(4.4)}
The group generated by\/ $a$ and\/ $b$ is isomorphic
to the restricted wreath product\/ $\mathbb Z\wr\mathbb Z$\textup;
more precisely\textup,
$$
\langle a,b\rangle=\langle b\rangle\wr\langle a\rangle,\quad
\langle b\rangle\cong\langle a\rangle\cong\mathbb Z.
$$
\end{lemma}

\begin{proof}
First of all, observe that
$\supp(a^{-k}ba^k)=(\alpha_k,\alpha_{k+1})$
for all $k\in\mathbb Z$
(see Lemma \ref{lemma:02.(2.2)}).
In particular, the supports of the maps $a^{-k}ba^k$,
$k\in\mathbb Z$, are pairwise disjoint,
and $\supp(a)=(0,r)$ is not equal to the support of any
element of the group $\langle\,a^{-k}ba^k\mid k\in\mathbb Z\,\rangle$.
Thus
$$
\langle\,a^{-k}ba^k\mid k\in\mathbb Z\,\rangle
=\bigoplus_{k\in\mathbb Z}\langle a^{-k}ba^k\rangle
\quad\text{and}\quad
\langle\,a^{-k}ba^k\mid k\in\mathbb Z\,\rangle\cap\langle a\rangle
=\{\id\},
$$
which implies
$\langle a,b\rangle
=\langle\,a^{-k}ba^{k}\mid k\in\mathbb Z\,\rangle\rtimes\langle a\rangle
=\langle b\rangle\wr\langle a\rangle$.
Since $\F$ has no torsion,
$\langle a\rangle\cong\langle b\rangle\cong\mathbb Z$.
\end{proof}

In the rest of this section,
$G$ is a subgroup of $\V(r,\mathbb R^*_+,\mathbb R)$
such that
$$
G\cap\F(r,\mathbb R^*_+,\mathbb R)=\F=\F(r,\Lambda,A)
$$
(it is even possible to generalize some results of this section
to the case when $G$ satisfies a weaker condition than this one).

\begin{proposition}
\label{proposition:02.(4.5)}
Let\/ $a$\textup, $b$\textup, and\/ $G$ be
the elements and the group defined above\textup.
Then there exists a first-order formula\/ $\phi$
with\/ $a$ and\/ $b$ as parameters and with exactly one free variable
such that\/ $\langle a,b\rangle$
is defined in\/ $G$ by\/ $\phi$\textup.
Moreover\textup, $\phi$ can be chosen based only on\/
$a$\textup, $b$\textup, and\/ $\F$
\textup(without knowing the whole of\/ $G$\textup{).}
\end{proposition}

In order to find such a $\phi$ and thus prove this proposition,
choose first of all $c,d\in \F$ and
$s,t\in\mathbb N$ such that
\begin{enumerate}
\item
	$c$ is a generator of the centralizer of $a$
	in $\F$,
\item
	$d$ is a generator of the centralizer of $b$
	in $\F_{(\alpha_0,\alpha_1)}$,
\item
	$a=c^s$ and $b=d^t$.
\end{enumerate}
The elements $c$ and $d$ and the numbers $s$ and $t$ exist by
Lemma \ref{lemma:11.(4.2)}.
By Lemma \ref{lemma:04.(2.4)}, $\supp(c)=(0,r)$
and $\supp(d)=(\alpha_0,\alpha_1)$.
Clearly $c,d\in\F^\uparrow$.

Observe that $a^{-k}da^k\in\F^\uparrow$ and
$\supp(a^{-k}da^k)=(\alpha_k,\alpha_{k+1})$ for all $k$.
Since for every $k$ the conjugation by $a^k$ is an automorphism
of $\F$ which sends
$\F_{(\alpha_0,\alpha_1)}$ onto
$\F_{(\alpha_{k},\alpha_{k+1})}$,
$a^{-k}da^k$ is a generator of the centralizer of $a^{-k}ba^k$
in $\F_{(\alpha_{k},\alpha_{k+1})}$
for every~$k$.

\begin{lemma}
\label{lemma:14.(4.6)}
The centralizer of\/ $\{\,a^{-k}ba^{k}\mid k\in\mathbb Z\,\}$
in\/ $\F$ is generated by\/
$\{\,a^{-k}da^{k}\mid k\in\mathbb Z\,\}$\textup.
\end{lemma}
\begin{proof}
The inclusion
$$
C_{\F}(\{\,a^{-k}ba^{k}\mid k\in\mathbb Z\,\})
\supset\langle\,a^{-k}da^{k}\mid k\in\mathbb Z\,\rangle
$$
is obvious, it remains show the inverse one.
So let $x$ be an arbitrary element of $\F$
which commutes with every $a^{-k}ba^{k}$, $k\in\mathbb Z$.

By Lemma \ref{lemma:02.(2.2)}, $x$ permutes each of the intervals
$(\alpha_k,\alpha_{k+1})$, $k\in\mathbb Z$.
By continuity and monotonicity, $(\alpha_k)x=\alpha_k$ for all $k$.
For every $k\in\mathbb Z$, let $y_k$ be the permutation of $[0,r)$
such that
$$
y_k|_{[\alpha_k,\alpha_{k+1}]}=x|_{[\alpha_k,\alpha_{k+1}]}
\quad\text{and}\quad
y_k|_{[0,\alpha_k]\cup[\alpha_{k+1},r)}
=\id_{[0,\alpha_k]\cup[\alpha_{k+1},r)}.
$$
Then for every $k\in\mathbb Z$,
$y_k\in\F_{(\alpha_{k},\alpha_{k+1})}$ and
$y_k$ commutes with $a^{-k}ba^{k}$.
Hence $y_k\in\langle a^{-k}da^{k}\rangle$ for all~$k$.

Choose $\beta,\gamma\in(0,r)$ such that
$x$ be affine on $[0,\beta]$ and on $[\gamma,r)$.
Then $x|_{[0,\beta]\cup[\gamma,r)}=\id_{[0,\beta]\cup[\gamma,r)}$.
Choose $n\in\mathbb N$ such that $\alpha_{-n}\in(0,\beta]$ and
$\alpha_{n+1}\in[\gamma,r)$.
Then $\supp(x)\subset(\alpha_{-n},\alpha_{n+1})$, and hence
$$
x=y_{-n}y_{-n+1}\dotsm y_{n-1}y_{n}
\in\langle\,a^{-k}da^{k}\mid k\in\mathbb Z\,\rangle.
$$
\end{proof}

\begin{lemma}
\label{lemma:15.(4.7)}
The centralizer of the element\/ $a$ and the centralizer of
the set\/ $\{\,a^{-k}ba^{k}\mid k\in\mathbb Z\,\}$
do not change when passing from\/ $\F$ to\/~$G$\textup.
\end{lemma}
\begin{proof}
By Lemmas \ref{lemma:08.(3.1)} and \ref{lemma:09.(3.2)},
all elements of $C_{G}(a)$ and of
$C_{G}(\{\,a^{-k}ba^{k}\mid k\in\mathbb Z\,\})$
are continuous.
Since an element of $G$ belongs to $\F$
if and only if it is continuous, the proof is complete.
\end{proof}

The group $\langle c\rangle$ is definable in $G$
with the parameter $a$
because it is the centralizer of $a$
(see Lemma \ref{lemma:15.(4.7)}).
The group $\langle a\rangle$ is definable in $G$
with the same parameter because
$$
\langle a\rangle=\bigl\{\,x^s\bigm|x\in\langle c\rangle\,\bigr\}.
$$
Thus the set $\{\,a^{-k}ba^{k}\mid k\in\mathbb Z\,\}$
is definable in $G$ with the parameters
$a$ and $b$, and hence so does its centralizer.
By Lemmas \ref{lemma:14.(4.6)} and \ref{lemma:15.(4.7)},
the centralizer of the set $\{\,a^{-k}ba^{k}\mid k\in\mathbb Z\,\}$
in $G$ is the group
$\langle\,a^{-k}da^{k}\mid k\in\mathbb Z\,\rangle$.
Since
$$
\langle\,a^{-k}ba^{k}\mid k\in\mathbb Z\,\rangle
=\bigl\{\,x^t\bigm|
x\in\langle\,a^{-k}da^{k}\mid k\in\mathbb Z\,\rangle\,\bigr\},
$$
the group $\langle\,a^{-k}ba^{k}\mid k\in\mathbb Z\,\rangle$ is
definable with the same parameters.
The group $\langle a,b\rangle$ is definable with
the parameters $a$ and $b$
since it is the semi-direct product of $\langle a\rangle$
and $\langle\,a^{-k}ba^{k}\mid k\in\mathbb Z\,\rangle$.

The following formula defines $\langle a,b\rangle$ in $G$
and depends only on $a$, $b$, $s$, and~$t$:
\begin{align*}
\phi(x)=\ulcorner\Bigl(\exists y,z\Bigr)
\Bigl(x=y^{s}z^{t}&\wedge ya=ay\\
&\wedge\Bigl(\forall w\Bigr)
\Bigl(wa=aw\rightarrow zw^{-s}bw^{s}=w^{-s}bw^{s}z\Bigr)\Bigr)\urcorner.
\end{align*}

We have proved Proposition \ref{proposition:02.(4.5)}.
We thus conclude:

\begin{corollary}
\label{corollary:01.(4.8)-proposition:02.(4.5)}
The group\/ $\F$ has subgroups isomorphic to\/ $\mathbb Z\wr\mathbb Z$
and definable with parameters in\/ $\F$\textup,
in\/ $\T$\textup, and in\/~$\V$\textup.
\end{corollary}


\section{Interpretations of the Arithmetic}
\label{section:interpretationarithmetic}

In this section we complete our proof of the interpretability of the
Arithmetic in $\F$, $\T$, and $\V$ with parameters.
Futhermore, in the case of the group $F$, or more generally,
of the group $\F$ with non-trivial cyclic $\Lambda$,
we exhibit an interpretation of the Arithmetic
which does not require parameters.

Apparently it is well known to specialists that every
finitely generated virtually solvable but not virtually abelian group
interprets the Arithmetic
(see \cite{Noskov:1983:etkpprg-rus,Noskov:1984:etfgasg-eng},
and also \cite{DelonSimonetta:1998:uwpspsf}).
For the reader's convenience, we present here our self-contained proof
for the group $\mathbb Z\wr\mathbb Z$.

\begin{lemma}
\label{lemma:16.(5.1)}
The group\/ $\mathbb Z\wr\mathbb Z$ interprets
the Arithmetic with parameters\textup.
More precisely\textup, if\/
$\mathbb Z\wr\mathbb Z=\langle b\rangle\wr\langle a\rangle$\textup,
$\langle b\rangle\cong\langle a\rangle\cong\mathbb Z$\textup,
then the bijection\/
$$
f\colon\{\,a^n\mid n\in\mathbb N\,\}\to\mathbb N,\quad a^n\mapsto n
$$
is an interpretation of\/
$(\mathbb N,+,\times)$ in\/ $(\mathbb Z\wr\mathbb Z,\times)$
with parameters\textup.
\end{lemma}
\begin{proof}
Denote
$$
G=\mathbb Z\wr\mathbb Z=\langle b\rangle\wr\langle a\rangle
\quad\text{et}\quad
H=\langle\,a^{-k}ba^{k}\mid k\in\mathbb Z\,\rangle.
$$
Recall that
$G=H\rtimes\langle a\rangle$,
and that $H$ is a free abelian group with the basis
$(a^{-k}ba^{k})_{k\in\mathbb Z}$.
It is easy to verify that $\langle a\rangle=C_G(a)$
and $H=C_G(b)$.

Consider the bijection
$$
g\colon \langle a\rangle\to\mathbb Z,\quad a^n\mapsto n.
$$
It will suffice to show that $g$ is an interpretation of
$(\mathbb Z,+,\times)$ in $(G,\times)$.
Indeed, $(\mathbb N,+,\times)$ is a substructure of
$(\mathbb Z,+,\times)$, and
$\mathbb N$ is $0$-definable
in $(\mathbb Z,+,\times)$ because,
by Lagrange's four-square theorem,
every positive integer is the sum of four squares.

The domain of $g$ is the centralizer of $a$,
hence definable.
The operation induced on $\langle a\rangle$
by the addition of $\mathbb Z$ via $g$ is
simply the restriction of the multiplication of $G$,
hence $0$-definable.
It remains to show that the operation induced on $\langle a\rangle$
by the multiplication of $\mathbb Z$ via $g$ is definable.

Observe the following facts:
\begin{enumerate}
\item
	for every $x\in G\setminus H$, $C_G(x)$ is cyclic,
\item
	for every $n\in\mathbb Z\setminus\{0\}$,
	$C_G(ba^n)=\langle ba^n\rangle$,
\item
	for every $n\in\mathbb Z$,
	$HC_G(ba^n)=H\langle a^n\rangle$.
\end{enumerate}
(The second fact is due to the homomorphism
$G\to\mathbb Z,\ a\mapsto0,\ b\mapsto1$.)

Denote by $|$ the relation of the divisibility
in $\mathbb Z$.
Observe that for all $m,n\in\mathbb Z$,
$$
m|n
\:\Leftrightarrow\:
HC_G(ba^m)\supset HC_G(ba^n).
$$
The relation $HC_G(bx)\supset HC_G(by)$ between $x,y\in G$
can be expressed by a first-order formula with the parameter $b$.
Hence the relation induced on $\langle a\rangle$ by $|$ via $g$
is definable.

The multiplication in $\mathbb Z$ is definable in terms of
the addition, the divisibility, and the constant $1$,
as can be seen from the following equivalences satisfied in~$\mathbb Z$:
\begin{align*}
n=k(k+1)&\leftrightarrow
\Bigl(\forall m\Bigr)\Bigl(n|m\leftrightarrow
k|m\wedge(k+1)|m\Bigr)\wedge(2k+1)|(2n-k),\\
n=kl&\leftrightarrow
(k+l)(k+l+1)=k(k+1)+l(l+1)+2n
\end{align*}
(see \cite[\S5a]{Robinson:1951:ur} for details).
Thus the operation induced on $\langle a\rangle$
by the multiplication of $\mathbb Z$
via $g$ is definable.
\end{proof}

Theorem A (see the Introduction)
is a corollary of Proposition \ref{proposition:02.(4.5)}
and Lemma \ref{lemma:16.(5.1)}.
In order to prove Theorem B, we shall construct
new $0$-interpretations of the Arithmetic in
groups of $\F$-kind:

\begin{proposition}
\label{proposition:03.(5.2)}
If\/ $\Lambda$ is cyclic\textup,
$\Lambda=\langle p\rangle$\textup,
then the map\/
$$
f\colon \{\,x\in\F\mid(0)x^{\prime+}=(r)x^{\prime-}>1\,\}\to\mathbb N,
\quad x\mapsto\log_p((0)x^{\prime+})
$$
is an interpretation of\/
$(\mathbb N,+,\times)$ in\/ $(\F,\times)$
without parameters\textup.
\end{proposition}

One of the main ideas of the proof of Proposition
\textup{\ref{proposition:03.(5.2)}} is the use of
the centralizers of pairs of elements.%
\footnote{Centralizers in $\F(r,\mathbb R^*_+,\mathbb R)$ have been
described by Brin and Squier \cite{BrinSquier:2001:pcrcgplhrl}.
Collin Bleak and others \cite{BleakGGHMNS:pp2007:} have
recently announced classification of all centralizers in
$\T(1,\langle n\rangle,\mathbb Z[\frac{1}{n}])$
and $\V(1,\langle n\rangle,\mathbb Z[\frac{1}{n}])$, $n=2,3,\dotsc$.}
The following lemma is similar to Theorem 5.5
in~\cite{BrinSquier:2001:pcrcgplhrl}:

\begin{lemma}
\label{lemma:17.(5.3)}
Let\/ $H$ be a subgroup of\/ $\F$\textup.
Then\/ $H$ is the centralizer of an element
if and only if\/ $H$ can be decomposed into a direct product
of subgroups\/ $H_1,\dotsc,H_n$\textup,
$n\in\mathbb N$\textup, such that there exist
$\alpha_0,\dotsc,\alpha_n\in A$ such that\/\textup:
\begin{enumerate}
\item
	$0=\alpha_0<\alpha_1<\dotsb<\alpha_n=r$\textup;
\item
	for every\/ $i=1,\dotsc,n$\textup,
	either\/ $H_i=\F_{(\alpha_{i-1},\alpha_i)}$\textup,
	or there exists\/ $x$ in\/ $\F_{(\alpha_{i-1},\alpha_i)}$ such that\/
	$H_i$ is the centralizer of $x$ in\/
	$\F_{(\alpha_{i-1},\alpha_i)}$\textup, and\/
	$H_i=\langle x\rangle$\textup;
\item
	for every\/ $i=1,\dotsc,n-1$\textup,
	if\/ $H_i=\F_{(\alpha_{i-1},\alpha_i)}$\textup, then\/
	$H_{i+1}\ne\F_{(\alpha_i,\alpha_{i+1})}$\textup.
\end{enumerate}
\end{lemma}
\begin{proof}
Let $x\in\F$ and $H=C_\F(x)$.
Choose $\alpha_0,\dotsc,\alpha_n$
such that $0=\alpha_0<\alpha_1<\dotsb<\alpha_n=r$ and
$$
\{\alpha_1,\dotsc,\alpha_{n-1}\}
=\{\,\alpha\in(0,r)\cap A\cap\fix(x)\mid
(\alpha)x^{\prime-}\ne1\text{ or }(\alpha)x^{\prime+}\ne1\,\}.
$$
For every $i=1,\dotsc,n$,
choose $x_i\in\F_{(\alpha_{i-1},\alpha_i)}$ such that
$x_i|_{(\alpha_{i-1},\alpha_i)}=x|_{(\alpha_{i-1},\alpha_i)}$,
and let $H_i$ be the centralizer of $x_i$
in $\F_{(\alpha_{i-1},\alpha_i)}$.
Note that $x=x_1\dotsm x_n$.
For every $i=1,\dotsc,n$, if $x_i\ne\id$, then $H_i$ is cyclic
(see Lemma \ref{lemma:11.(4.2)}).

Similarly to Lemma \ref{lemma:03.(2.3)}, it is easy to prove that
every element $y$ of $H$ permutes
the set $\{\alpha_0,\dotsc,\linebreak[0]\alpha_{n-1}\}$,
and hence, as this set is finite, $y$ fixes all its elements.
Thus $H=H_1\times\dotsm\times H_n$.

Conversely, suppose that $H=H_1\times\dotsm\times H_n$,
$n\in\mathbb N$, and that the subgroups $H_1,\dotsc,H_n$
and the points $\alpha_0,\dotsc,\alpha_n\in A$
are as in the statement of this lemma.
For every $i=1,\dotsc,n$, choose $x_i\in\F_{(\alpha_{i-1},\alpha_i)}$
such that $H_i$ be the centralizer of $x_i$
in $\F_{(\alpha_{i-1},\alpha_i)}$.
Then $H=C_\F(x_1\dotsm x_n)$.
\end{proof}

\begin{lemma}
\label{lemma:18.(5.4)}
Let\/ $H$ be a subgroup of\/ $\F$\textup.
Then\/ $H$ is the centralizer of a pair of elements\/
\textup(possibly equal\/\textup)
if and only if\/ $H$ can be decomposed into a direct product
of subgroups\/ $H_1,\dotsc,H_n$\textup,
$n\in\mathbb N$\textup, such that there exist
$\alpha_0,\dotsc,\alpha_n\in A$ such that\/\textup:
\begin{enumerate}
\item
	$0=\alpha_0<\alpha_1<\dotsb<\alpha_n=r$\textup;
\item
	for every\/ $i=1,\dotsc,n$\textup,
	either\/ $H_i=\{\id\}$\textup,
	or\/ $H_i=\F_{(\alpha_{i-1},\alpha_i)}$\textup,
	or there exists\/ $x$ in\/ $\F_{(\alpha_{i-1},\alpha_i)}$ such that\/
	$H_i$ is the centralizer of $x$ in\/
	$\F_{(\alpha_{i-1},\alpha_i)}$\textup, and\/
	$H_i=\langle x\rangle$\textup;
\item
	for every\/ $i=1,\dotsc,n-1$\textup,
	if\/ $H_i=\F_{(\alpha_{i-1},\alpha_i)}$\textup, then\/
	$H_{i+1}\ne\F_{(\alpha_i,\alpha_{i+1})}$\textup.
\end{enumerate}
\end{lemma}
\begin{proof}
This lemma is an easy corollary of Lemma \ref{lemma:17.(5.3)}
and the fact that for all $\alpha,\beta\in A$
such that $0<\alpha<\beta<r$, there exist $x,y\in\F$ such that
$\supp(x)=\supp(y)=(\alpha,\beta)$ and $xy\ne yx$,
and hence the centralizer of $\{x,y\}$ in $\F_{(\alpha,\beta)}$
is trivial (see Lemmas \ref{lemma:11.(4.2)} and \ref{lemma:12.(4.3)}).
\end{proof}

\begin{proof}
[Proof of Proposition\/ \textup{\ref{proposition:03.(5.2)}}]
Without loss of generality, suppose that $p>1$.

Denote by $\F^\circ$ the subgroup of $\F$ formed by the elements
that are identities in neighborhoods of $0$ and $r$:
$$
\F^\circ=\{\,x\in\F\mid
(0)x^{\prime+}=(r)x^{\prime-}=1\,\}.
$$
Denote by $B$ the domain of $f$:
$$
B=\{\,x\in\F\mid(0)x^{\prime+}=(r)x^{\prime-}>1\,\}.
$$
Note that for any two elements $x$ and $y$ of $B$,
$f(x)=f(y)$ if and only if $xy^{-1}\in\F^\circ$
(see Lemma~\ref{lemma:10.(4.1)}).

Lemma \ref{lemma:12.(4.3)} allows to conclude that
$f$ is a surjection onto $\mathbb N$
(see Figure \ref{figure:3} for example).
\begin{figure}
\includegraphics{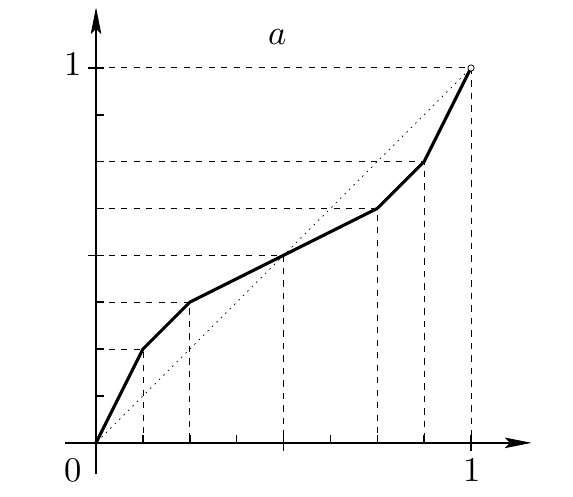}
\caption{An element $a$ of $F$ such that $(0)a^{\prime+}=(1)a^{\prime-}=2$.}
\label{figure:3}
\end{figure}
It will suffice to show that the set $B$, the group $\F^\circ$,
and the relations induced on $B$ via $f$ by
the addition and the divisibility of $\mathbb N$
are all $0$-definable (see Remark \ref{remark:04.(2.10)}).
Indeed, the multiplication is $0$-definable in $\mathbb N$
in terms of the addition and the divisibility---see
\cite[\S4b]{Robinson:1951:ur} or the proof of Lemma~\ref{lemma:16.(5.1)}.

Define
\begin{align*}
S&=\{\,\F_{(\alpha,\beta)}\mid\alpha,\beta\in A,\
0\le\alpha<\beta\le r\,\},\\
S_0&=\{\,\F_{(\alpha,\beta)}\mid\alpha,\beta\in A,\
0<\alpha<\beta<r\,\},\\
S_1&=\{\,\F_{(0,\beta)}\mid\beta\in A,\ 0<\beta<r\,\}\\
&\qquad\qquad\cup
\{\,\F_{(\alpha,r)}\mid\alpha\in A,\ 0<\alpha<r\,\}.
\end{align*}
Then not only is the family $S$ \emph{uniformly definable},
but also
there exist two first-order formulae
$\phi(x_1,x_2,x_3)$ and $\psi(x_1,x_2)$
(in the language of groups) \emph{without parameters\/} such that:
\begin{enumerate}
\item
	for all $\alpha,\beta\in A$ such that
	$0\le\alpha<\beta\le r$, there exist $x,y\in\F$ such that
	$\F\models\psi(x,y)$ and
	$\F_{(\alpha,\beta)}=\{\,z\in\F\mid\F\models\phi(x,y,z)\,\}$;
\item
	for all $x,y\in\F$ such that
	$\F\models\psi(x,y)$, there exist $\alpha,\beta\in A$ such that
	$0\le\alpha<\beta\le r$ and
	$\F_{(\alpha,\beta)}=\{\,z\in\F\mid\F\models\phi(x,y,z)\,\}$.
\end{enumerate}
Indeed, a subset of $\F$ is of the form
$\F_{(\alpha,\beta)}$, where $\alpha,\beta\in A$ and
$0\le\alpha<\beta\le r$,
if and only if it is the centralizer of a pair,
is not abelian, and cannot be decomposed as the direct product
of two other centralizers of pairs
(see Lemma \ref{lemma:18.(5.4)}).
All this can be expressed in the first-order language.
The families $S_0$ and $S_1$ of subsets of $\F$ are
``uniformly definable without parameters'' in the same sense as $S$
because
\begin{align*}
S_0&=\{\,H_1\cap H_2\mid
H_1,H_2\in S,\ H_1\not\subset H_2,\text{ and }H_2\not\subset H_1\,\},\\
S_1&=S\setminus(S_0\cup\{\F\}).
\end{align*}
Since $\F^\circ=\bigcup S_0$, this subgroup is $0$-definable.%
\footnote{In the case when $\F=F$,
there exists a more natural proof that 
$\F^\circ$ is definable in $\F$
since in this case $\F^\circ=F^\circ=[F,F]=[\F,\F]$
(see \cite[Theorem 4.1]{CannonFloydParry:1996:inRTg}),
and every element of $[\F,\F]$ is the product of two commutators
(see the Appendix).
In general $\F^\circ$ and $[\F,\F]$ are not equal
(see \cite[Section 4D]{Brown:1987:fpg}).}

Define
\begin{align*}
E&=\{\,x\in\F\mid(0)x^{\prime+}\ne1\text{ or }(r)x^{\prime-}\ne1\,\}
=\F\setminus\bigcup S_0=\F\setminus\F^\circ,\\
E_2&=\{\,x\in\F\mid(0)x^{\prime+}\ne1\text{ and }(r)x^{\prime-}\ne1\,\}
=\F\setminus\bigcup S_1=E\setminus\bigcup S_1.
\end{align*}
These sets are $0$-definable.

Define
\begin{align*}
P^+&=\{\,x\in\F\mid(0)x^{\prime+}>1\text{ and }(r)x^{\prime-}>1\,\},\\
P^-&=\{\,x\in\F\mid(0)x^{\prime+}<1\text{ and }(r)x^{\prime-}<1\,\},
\end{align*}
and $P=P^+\cup P^-$.
These sets are $0$-definable:
for all $x\in\F$,
$$
x\in P^+
\:\Leftrightarrow\:
\Bigl(\exists X\in S_0\Bigr)\,\Bigl(\forall Y\in S_1\Bigr)\,
\Bigl(Y\supset X\rightarrow x^{-1}Yx\subsetneqq Y\Bigr),
$$
and
$$
x\in P^-
\:\Leftrightarrow\:
\Bigl(\exists X\in S_0\Bigr)\,\Bigl(\forall Y\in S_1\Bigr)\,
\Bigl(Y\supset X\rightarrow x^{-1}Yx\supsetneqq Y\Bigr).
$$
It has been used here that for all $\alpha,\beta$
and for all $x\in\F$,
$$
x^{-1}\F_{(\alpha,\beta)}x=\F_{((\alpha)x,(\beta)x)}.
$$

Define
$$
U=\{\,x\in\F\mid(0)x^{\prime+}=((r)x^{\prime-})^{-1}\in\{p^{\pm1}\}\,\}.
$$
Then $U$ is $0$-definable: $U\subset E_2\setminus P$, and
for all $x\in E_2\setminus P$,
\begin{multline*}
x\in U
\:\Leftrightarrow\:
\Bigl(\forall y\in P\Bigr)\,
\Bigl(\forall z\in\F^\circ\Bigr)\,
\Bigl(\exists w_1,w_2\in C_\F(xz)\Bigr)\\
\Bigl((w_1w_2^{-1}\in E_2)\wedge (yw_1,yw_2\notin E_2)\Bigr).
\end{multline*}
To facilitate reading of the last formula, we remark that
$xy^{-1}\notin E_2$ means exactly that either
$(0)x^{\prime+}=(0)y^{\prime+}$, or $(r)x^{\prime-}=(r)y^{\prime-}$.
Then the implication ``$\Rightarrow$'' is easy to check:
such $w_1,w_2$ can always be found even in $\langle xz\rangle$.
The direction ``$\Leftarrow$'' is less obvious, we shall rather
prove its contrapositive.

Let $x\in(E_2\setminus P)\setminus U$.
Without loss of generality, suppose that
$(0)x^{\prime+}>p$, and hence $(r)x^{\prime-}<1$.
Choose $y\in P$ such that $(0)y^{\prime+}=p$
(see Lemma \ref{lemma:12.(4.3)}).
Let $\gamma\in(0,r)\cap A$, and choose $x_1,x_2\in\F$ such that:
\begin{enumerate}
\item
	$\supp(x_1)=(0,\gamma)$, $\supp(x_2)=(\gamma,r)$,
\item
	$(0)x_1^{\prime+}=(0)x^{\prime+}$,
	$(\gamma)x_1^{\prime-}=p^{-1}$,
\item
	$(r)x_2^{\prime-}=(r)x^{\prime-}$,
	$(\gamma)x_2^{\prime+}=p$.
\end{enumerate}
Then the centralizer of $x_1x_2$ is the direct product
$\langle x_1\rangle\times\langle x_2\rangle$
(see Lemmas \ref{lemma:03.(2.3)} and \ref{lemma:10.(4.1)}).
Choose $z\in\F^\circ$ such that $xz=x_1x_2$.
Now suppose that
$$
w_1,w_2\in C_\F(xz)=\langle x_1\rangle\times\langle x_2\rangle.
$$
Then $(0)(yw_1)^{\prime+}\ne1$ and $(0)(yw_2)^{\prime+}\ne1$.
Suppose that $yw_1,yw_2\notin E_2$.
Then $(r)(yw_1)^{\prime-}=(r)(yw_2)^{\prime-}=1$,
and hence $(r)(w_1w_2^{-1})^{\prime-}=1$,
or in other terms $w_1w_2^{-1}\notin E_2$.
Thus $U$ is indeed $0$-definable.

For all $x\in P$,
$$
(0)x^{\prime+}=(r)x^{\prime-}
\:\Leftrightarrow\:
\Bigl(\forall y\in U\Bigr)\,
\Bigl(\exists z\in C_\F(y)\Bigr)\,
\Bigl(xz,xz^{-1}\notin E_2\Bigr).
$$
Thus, the set $B$ is $0$-definable.

The $f$-preimage of the graph of the addition of $\mathbb N$
is $0$-definable since it is the graph of the multiplication
modulo $\F^\circ$: for all $x,y,z\in B$,
$$
(x)f+(y)f=(z)f
\:\Leftrightarrow\:
xyz^{-1}\in\F^\circ.
$$

It remains to show that the $f$-preimage of the graph of the
divisibility of $\mathbb N$ is $0$-definable.
This is indeed the case:
for all $x,y\in B$,
$$
(x)f|(y)f
\:\Leftrightarrow\:
\Bigl(\forall z\in\F^\circ\Bigr)\,
\Bigl(\exists w\in C_\F(xz)\Bigr)\,
\Bigl(yw\in\F^\circ\Bigr).
$$
The implication ``$\Leftarrow$'' here is the least evident.
In order to establish its contrapositive, one can take
$\gamma\in(0,r)\cap A$, choose $x_1,x_2\in\F$ such that:
\begin{enumerate}
\item
	$\supp(x_1)=(0,\gamma)$, $\supp(x_2)=(\gamma,r)$,
\item
	$(0)x_1^{\prime+}=(0)x^{\prime+}$,
	$(\gamma)x_1^{\prime-}=p^{\pm1}$,
\item
	$(r)x_2^{\prime-}=(r)x^{\prime-}$,
	$(\gamma)x_2^{\prime+}=p^{\pm1}$,
\end{enumerate}
and choose $z\in\F^\circ$ such that $xz=x_1x_2$
(and hence $C_\F(xz)=\langle x_1\rangle\times\langle x_2\rangle$).
\end{proof}

Theorem B is a corollary of
Proposition \textup{\ref{proposition:03.(5.2)}}.


\section{Undecidability}
\label{section:undecidability}

In this section we deduce from
\cite[Theorem 9]{MostowskiRobTar:1971:ueua}
(see Theorem \ref{theorem:01.(6.3)} below)
that the elementary theory of every structure of finite signature
that interprets the Arithmetic with parameters is
hereditarily undecidable.

In what follows, the constants are treated as functions of
\emph{arity\/} $0$,
and similarly constant symbols are viewed as a particular case
of function symbols.
If $\sigma$ is a relation symbol, function symbol, or a constant symbol,
its arity shall be denoted by $\arity(\sigma)$.

We shall say that an $n$-ary relation $R$ on a set $B$ is
\emph{compatible\/} with an equivalence relation $E$ on $B$ if
$R$ is induced by a relation (of the same arity) on $B/E$,
\ie\ if belonging of an $n$-tuple $(b_1,\dotsc,b_n)$ to $R$
is completely determined by the classes of $E$-equivalence of
$b_1,\dotsc,b_n$.

Let $\Sigma$ and $\Gamma$ be two signatures, and let
$N$ be a $\Gamma$-structure.
Let $n\in\mathbb N$, $B$ be a definable subset of $N^{n}$,
and $E$ be an equivalence relation on $B$, also definable in $N$.
Let $\phi$, $\psi$, and $\xi_\sigma$ for all $\sigma\in\Sigma$
be $\Gamma$-formulae with parameters from $N$ such that:
\begin{enumerate}
\item
	$\phi$ defines $B$,
\item
	$\psi$ defines $E$,
\item
	for every relation symbol $\sigma\in\Sigma$,
	the formula $\xi_\sigma$
	defines a relation on $B$ compatible with $E$
	of arity $\arity(\sigma)n$, and
\item
	for every function or constant symbol $\sigma\in\Sigma$,
	the formula $\xi_\sigma$
	defines a relation on $B$
	of arity $(\arity(\sigma)+1)n$ which is compatible with $E$ and
	which defines the graph of an operation on $B/E$
	of arity $\arity(\sigma)$.
\end{enumerate}
Then denote by $\Int{\Sigma}{N,\phi,\psi,(\xi_\sigma)_{\sigma\in\Sigma}}$
the $\Sigma$-structure naturally defined on $B/E$ by the family
$(\xi_\sigma)_{\sigma\in\Sigma}$.

\begin{remark}
\label{remark:07.(6.1)}
In the same notation, the natural projection
$p\colon B\to B/E$ is an interpretation
of $\Int{\Sigma}{N,\phi,\psi,(\xi_\sigma)_{\sigma\in\Sigma}}$ in~$N$.
\end{remark}

\begin{lemma}
\label{lemma:19.(6.2)}
Let\/ $M$ and\/ $N$ be two structures of finite signatures such that\/
$\Th(M)$ is hereditairily undecidable\textup, and\/ $N$ interprets\/
$M$ with parameters\textup.
Then\/ $\Th(N)$ is hereditarily undecidable as well\textup.
\end{lemma}
\begin{proof}
Denote the signature of $M$ by $\Sigma$
and the signature of $N$ by~$\Gamma$.

It suffices to consider only the case when $\Sigma$ contains
no function or constant symbols.
Indeed, let $\Sigma'$ be the signature obtained from $\Sigma$
by replacing every $n$-ary function symbol $f$ ($n\ge0$)
by an $(n+1)$-ary relation symbol $f'$.
For every $\Sigma$-structure $M$, denote by $M'$ the $\Sigma'$-structure
on the underlying set of $M$ in which every new relation symbol
of $\Sigma'$ is interpreted by the graph of the function in $M$
named by the corresponding symbol of $\Sigma$,
and all the other symbols of $\Sigma'$ are interpreted in $M'$
exactly like in $M$.
For every $\Sigma$-theory $S$,
denote by $S'$ the $\Sigma'$-theory
of the class $\{\,M'\mid M\models S\,\}$.
It is easy to see that:
\begin{enumerate}
\item
	the (class-)map $M\mapsto M'$,
	where $M$ is a $\Sigma$-model of $S$,
	is a (class-) bijection between
	$\Mod_{\Sigma}S$ and $\Mod_{\Sigma'}S'$ for every
	$\Sigma$-theory $S$,
\item
	if $O_{\Sigma}$ denotes the minimal $\Sigma$-theory, then
	the map $S\mapsto S'$
	of $\Sigma$-theories to $\Sigma'$-theories
	is a bijection between all the $\Sigma$-theories
	and all the $\Sigma'$-theories containing $O_{\Sigma}'$
	($O_{\Sigma}'$ expresses simply that the new relation symbols
	of $\Sigma'$ are to be interpreted by graphs of functions),
\item
	for every $\Sigma$-structure $M$,
	a set is $0$-definable in $M$ if and only if it is $0$-definable in $M'$,
	or in other words
	$\id_M$ is a $0$-interpretation of $M$ in $M'$ and of $M'$ in~$M$.
\end{enumerate}
It is easy to provide an algorithm which converts every
$\Sigma$-sentence $\phi$ into a $\Sigma'$-sentence
$\psi$ such that for every $\Sigma$-structure $M$,
$M\models\phi$ if and only if $M'\models\psi$,
and it is equally easy to provide and algorithm which converts every
$\Sigma'$-sentence into a $\Sigma$-sentence equivalent in the same sense.
Therefore, for every $\Sigma$-theory $S$,
$S'$ is essentially undecidable if and only if such is $S$.
Hence we suppose without loss of generality that
$\Sigma$ contains only relation symbols.

Let $(n,f)$ be an interpretation of $M$ in $N$.
Let $a_1,\dotsc,a_m$ be a sequence of parameters from $N$
sufficient to define the domain and the kernel of $f$ and
the $f$-preimage of the graph of every relation of $M$
($\Sigma$ is finite);
denote $\bar a=(a_1,\dotsc,a_m)$.
Let $x_1,\dotsc,x_m$, $y_1,\dotsc,y_n$,
$y_{11},\dotsc,y_{1n}$, $y_{21},\dotsc,y_{2n}$, \dots\
be distinct variables, and denote
$\bar x=(x_1,\dotsc,x_m)$,
$\bar y=(y_1,\dotsc,y_n)$,
$\bar y_1=(y_{11},\dotsc,y_{1n})$, and so on.
Let $\phi=\phi(\bar x,\bar y)$,
$\psi=\psi(\bar x,\bar y_1,\bar y_2)$,
and $\xi_\sigma=\xi_\sigma(\bar x,\bar y_1,\dotsc,\bar y_k)$
for every $\sigma\in\Sigma$ of arity $k$
be $\Gamma$-formulae such that:
\begin{enumerate}
\item
	$\phi(\bar a,\bar y)$ defines the domain of $f$
	(which is a subset of $N^n$),
\item
	$\psi(\bar a,\bar y_1,\bar y_2)$ defines the kernel of $f$
	(which is an equivalence relation on the domain),
\item
	for every symbol $\sigma\in\Sigma$,
	the formula $\xi_\sigma(\bar a,\bar y_1,\dotsc,\bar y_{\arity(\sigma)})$
	defines the $f$-preimage of the graph of the relation of $M$ named
	by~$\sigma$.
\end{enumerate}
Observe that the bijection induced by $f$ between the quotient of its domain
by its kernel and its image is an isomorphism
$$
\Int{\Sigma}{N,\phi(\bar a,\bar y),
\psi(\bar a,\bar y_1,\bar y_2),
(\xi_\sigma(\bar a,\bar y_1,\dotsc,
\bar y_{\arity(\sigma)}))_{\sigma\in\Sigma}}\stackrel{\cong}{\to}M.
$$

In what follows, let $\bar c=(c_1,\dotsc,c_m)$ be a
sequence of new constant symbols.
Write $(\Gamma,\bar c)$ to denote the signature obtained
from $\Gamma$ by adding $c_1,\dotsc,c_m$
(as constant symbols).

Let $\tau=\tau(\bar x)$ be a $\Gamma$-formula such that
the $(\Gamma,\bar c )$-sentence $\tau(\bar c)$ expresses that:
\begin{enumerate}
\item
	$\psi(\bar c,\bar y_1,\bar y_2)$
	defines an equivalence relation on
	the set defined by $\phi(\bar c,\bar y)$,
\item
	for every relation symbol $\sigma\in\Sigma$,
	the formula
	$\xi_\sigma(\bar c,\bar y_1,\dotsc,\bar y_{\arity(\sigma)})$
	defines a relation
	on the set defined by $\phi(\bar c,\bar y)$
	which is compatible with the equivalence relation
	defined by $\psi(\bar c,\bar y_1,\bar y_2)$.
\end{enumerate}
Clearly all this can be expressed in the first-order language.%
\footnote{One can take as $\tau(\bar c)$
the conjunction of the \emph{admissibility conditions\/} in the sense of
\cite[Section 5.3]{Hodges:1993:mt}.}
Note that
$$
N\models\tau(\bar a).
$$

Choose a recursive (\ie\ computable by an algorithm) map $t$
from the set of $\Sigma$-sentences
to the set of $\Gamma$-formulae all of whose free variables are among
$x_1,\dotsc,x_m$ such that
for every $\Sigma$-sentence $\alpha$ and
every $(\Gamma,\bar c)$-structure $L$ such that $L\models\tau(\bar c)$,
$$
\Bigl(L\models\alpha^t(\bar c)\Bigr)
\Leftrightarrow
\Bigl(\Int{\Sigma}{L,\phi(\bar c,\bar y),\psi(\bar c,\bar y_1,\bar y_2),
(\xi_\sigma(\bar c,\bar y_1,\dotsc,
\bar y_{\arity(\sigma)}))_{\sigma\in\Sigma}}\models\alpha\Bigr).
$$
It is easy to construct such a $t$ that uses the formula
$\phi$ to relativize the quantifiers,
the formula $\psi$ to replace $\ulcorner=\urcorner$, and the formula
$\xi_\sigma$ to replace each $\sigma\in\Sigma$.%
\footnote{In \cite[Section 5.3]{Hodges:1993:mt}
such $t$ is called a \emph{reduction map}.}
Here is an example, where $\sigma$ is a binary relation symbol:
\begin{align*}
\alpha&=\ulcorner\Bigl(\forall y_1,y_2,y_3\Bigr)\,
\Bigl(\sigma(y_1,y_2)\wedge\sigma(y_1,y_3)
\rightarrow y_2=y_3\Bigr)\urcorner,\\
\alpha^t(\bar x)&=\ulcorner\Bigl(\forall\bar y_1,\bar y_2,\bar y_3\Bigr)\,
\Bigl(\phi(\bar x,\bar y_1)\wedge\phi(\bar x,\bar y_2)
\wedge\phi(\bar x,\bar y_3)\\
&\qquad\qquad\rightarrow
\Bigl(\xi_\sigma(\bar x,\bar y_1,\bar y_2)
\wedge\xi_\sigma(\bar x,\bar y_1,\bar y_3)
\rightarrow\psi(\bar x,\bar y_2,\bar y_3)\Bigr)\Bigr)\urcorner.
\end{align*}
(As is customary, we do not show all the parentheses;
they should be added according to the standard rules.)

Note the following properties of $t$:
\begin{enumerate}
\item
	for every $\Sigma$-sentence $\alpha$,
	$$
	\Bigl(M\models\alpha\Bigr)\Leftrightarrow
	\Bigl(N\models\alpha^t(\bar a)\Bigr);
	$$
\item
	for every $\Sigma$-sentence $\alpha$,
	$$
	\Bigl(\vdash_\Sigma\alpha\Bigr)\Rightarrow
	\Bigl(\tau(\bar c)\vdash_{(\Gamma,\bar c)}\alpha^t(c)\Bigr);
	$$
\item
	for every $\Sigma$-sentences $\alpha$ and $\beta$,
	\begin{align*}
	\tau(\bar c)&\vdash_{(\Gamma,\bar c)}(\alpha\wedge\beta)^t(c)
	\leftrightarrow\alpha^t(c)\wedge\beta^t(c),\\
	\tau(\bar c)&\vdash_{(\Gamma,\bar c)}(\neg\alpha)^t(c)
	\leftrightarrow\neg\alpha^t(c),
	\end{align*}
	and the same for the other boolean operations;
\item
	for every $(\Gamma,\bar c)$-theory $T$ such that
	$T\vdash_{(\Gamma,\bar c)}\tau(\bar c)$, the set
	$$
	\{\,\text{$\Sigma$-sentence }\alpha
	\mid T\vdash_{(\Gamma,\bar c)}\alpha^t(\bar c)\,\}
	$$
	is a $\Sigma$-theory.
\end{enumerate}

Suppose now that $\Th(N)$ were not hereditarily undecidable.
Then let $T$ be a decidable $\Gamma$-subtheory of $\Th(N)$.
Let $S$ be the set of all the $\Gamma$-formulae $\alpha(\bar x)$
such that
$$
T\vdash_\Gamma\Bigl(\forall\bar x\Bigr)\,
\Bigl(\tau(\bar x)\rightarrow\alpha(\bar x)\Bigr).
$$
Then $\tau\in S$, $N\models\alpha(\bar a)$ for all $\alpha\in S$,
$S$ is a recursive (decidable) set,
and $\{\,\alpha(\bar c)\mid\alpha(\bar x)\in S\,\}$
is a $(\Gamma,\bar c)$-theory.
Let $U$ be the preimage of $S$ under $t$.
Then $U\subset\Th(M)$, $U$ is a $\Sigma$-theory,
and $U$ is decidable in contradiction with the hereditary undecidability
of~$\Th(M)$.
\end{proof}

\begin{theorem}[Mostowski, Tarski, \cite{MostowskiTarski:1949:uaitr}]
\label{theorem:01.(6.3)}
The elementary theory of the Arithmetic\/ $(\mathbb N,+,\times)$
has an essentially undecidable finitely axiomatized subtheory\textup.
\end{theorem}

For an improved proof of this fact,
see \cite[Theorem~9]{MostowskiRobTar:1971:ueua}.

\begin{proposition}
\label{proposition:04.(6.4)}
Let\/ $M$ be a structure of finite signature
which interprets the Arithmetic\/ $(\mathbb N,+,\times)$
with parameters\textup.
Then\/ $\Th(M)$ is hereditarily undecidable\textup.
\end{proposition}
\begin{proof}
This is a corollary of Theorem \ref{theorem:01.(6.3)}
and Lemmas \ref{lemma:07.(2.13)} and~\ref{lemma:19.(6.2)}.
\end{proof}

Theorem C (see the Introduction)
follows now from Theorem A and Proposition \ref{proposition:04.(6.4)}.


\section{Open questions}
\label{section:questions}

We conclude with two questions which, to our knowledge, are open.

Thompson's group $F$ is definable in $T$.%
\footnote{We shall not show this here.}
However, it is not known to the authors whether $F$ is definable in~$V$.

\begin{question}
\label{question:1}
Is Thompson's group $F$ definable in $V$ with parameters?
\end{question}

We have already shown that the Arithmetic is interpreted in $F$.
In addition, $F$ is interpreted in the Arithmetic because
the word problem for $F$ is decidable.
(Matiyasevich's theorem, see
\cite{Matijasevich:1970:dpm-rus}
or \cite[Theorem 8.1]{Davis:1973:htpiu},
implies that every recursive or recursively enumerable subset of
$\mathbb N^n$, $n\in\mathbb N$, is definable in the Arithmetic;
thus every reasonable encoding of elements of $F$
by positive integers gives an interpretation of $F$ in the Arithmetic.)
However, even if a structure interprets the Arithmetic and, reciprocally,
is interpreted in the Arithmetic, they
are not necessarily \emph{bi-interpretable\/}
(see \cite[Theorem 6]{Khelif:2007:bisQFAegrac-fr} or
\cite[Theorem 7.16]{Nies:2007:dg}), hence the question:

\begin{question}
\label{question:2}
Is the group $F$ bi-in\-ter\-pret\-a\-ble with the Arithmetic
with parameters?
\end{question}

We say that a structure $S$ is
\emph{categorically finitely axiomatized\/}
in a class $C$ of structures of the same signature
if $S\in C$ and there exists a first-order sentence $\phi$
such that $S\models\phi$ and every structure in $C$ that satisfies $\phi$
is isomorphic to $S$.%
\footnote{In the case when the class $C$ consists of all the
finitely generated structures of a certain class,
André Nies \cite{Nies:2007:dg} and
Anatole Khélif \cite{Khelif:2007:bisQFAegrac-fr}
used the term ``\emph{quasi finitely axiomatized\/}''
in more or less the same sense as we use
``categorically finitely axiomatized.''
However, the definitions in \cite{Nies:2007:dg} and
\cite{Khelif:2007:bisQFAegrac-fr}
are not precise because
the property of finite generation is not intrinsic and depends
on the class in which a given structure is considered.}
According to Anatole Khélif \cite{Khelif:2007:bisQFAegrac-fr},
bi-interpretability with the Arithmetic can be used 
to demonstrate categoric finite axiomatization in classes
of finitely generated structures of finite signature.
Thomas Scanlon \cite{Scanlon:2008:ifgfbN}
has recently established bi-interpretability
of the Arithmetic with all finitely generated fields
and used it to show that all such fields are
categorically finitely axiomatized
within the class of finitely generated fields, and thus
Pop's conjecture holds true:
\begin{quote}
two finitely generated fields are elementary equivalent
if and only if they are isomorphic.
\end{quote}
André Nies raised the question
whether there exists a finitely generated simple group
categorically finitely axiomatized among all the
finitely generated simple groups
\cite[Question 7.8]{Nies:2007:dg}.


\section{Appendix}
\label{section:appendix}

Here we show that every element of the derived subgroup of $\F$
is the product of two commutators, and hence
$[\F,\F]$ is $0$-definable in $\F$.
The proof of this fact, which, incidentally,
has not been used in this paper, is known to specialists.
However, this fact is closely related to the definable structure of
the groups that we study here,
and apparently it does not appear anywhere else in the literature.
For Thompson's group $F$, this result is probably part of folklore;
we learned its proof from Matthew Brin, who had slightly modified the argument
of Keith Dennis and Leonid Vaserstein
\cite[Proposition 1(c)]{DennisVaserstein:1989:clg}.
Our argument is just a trivial generalization.

As in the proof of Proposition \ref{proposition:03.(5.2)},
define
$$
\F^\circ=\{\,x\in\F\mid
(0)x^{\prime+}=(r)x^{\prime-}=1\,\}.
$$

\begin{proposition}
\label{proposition:05.(8.1)}
Every element of\/ $[\F,\F]$ is the product of two commutators
in\/ $\F$\textup, and even in\/ $\F^\circ$\/\textup:
$$
[\F,\F]=\{\,[x_1,x_2][x_3,x_4]\mid x_1,x_2,x_3,x_4\in\F^\circ\,\}.
$$
\end{proposition}
\begin{proof}
First of all recall that $[\F,\F]\subset\F^\circ$
(see Lemma \ref{lemma:10.(4.1)}).

Consider any two elements $x,y\in\F$ and their commutator $c=[x,y]$.
Choose $\alpha_1,\alpha_2,\beta_1,\beta_2\in A$ such that
$0<\alpha_2<\alpha_1<\beta_1<\beta_2<r$ and
$\supp(c)\subset(\alpha_1,\beta_1)$.
Then there exists an endomorphism $h\colon\F\to\F_{(\alpha_2,\beta_2)}$
such that $h$ is the identity on $\F_{(\alpha_1,\beta_1)}$.
Indeed, let $s\colon[0,r)\to[0,r)$ be a map which is the
identity on $[\alpha_1,\beta_1]$ and which is affine on
$[0,\alpha_1]$ and on $[\beta_1,r)$
with the slope $p\ll1$, $p\in\Lambda$;
then it is possible to take as $h$ the conjugation by $s$
composed with the natural embedding of permutations of the interval
$[(0)s,(r)s)$ into permutation of $[0,r)$
(this $h$ is even injective).
If $h$ is such an endomorphism, if $(x)h=x'$, and if $(y)h=y'$, then
$$
c=(c)h=[x',y']\quad\text{and}\quad
\supp(x')\cup\supp(y')\subset(\alpha_2,\beta_2).
$$
The above has two consequences of interest for us:
\begin{enumerate}
\item
    if $x,y\in\F$, then there exist $x',y'\in\F^\circ$
    such that $[x,y]=[x',y']$;
\item
    if $c_1,c_2,\dotsc,c_n$ are commutators in $\F$,
    and $I_1,I_2,\dotsc,I_n$ are pairwise disjoint closed subintervals
    of $[0,r)$ such that
    $\supp(c_i)\subset I_i$ for all $i=1,\dotsc,n$,
    then the product $c_1c_2\dotsm c_n$ is also a commutator.
\end{enumerate}

Now let $c_1$, $c_2$, and $c_3$ be three arbitrary commutators in $\F$.
Choose $\alpha,\beta\in A\cap(0,r)$ such that
$$
\supp(c_1)\cup\supp(c_2)\cup\supp(c_3)\subset(\alpha,\beta).
$$
Choose $b\in\F$ such that $(\alpha)b>\beta$.
Then $0<(\beta)b^{-1}<\alpha<\beta<(\alpha)b<r$.
Since
$$
\supp(c_2^b)\subset((\alpha)b,r)
\quad\text{and}\quad
\supp(c_3^{b^{-1}})\subset(0,(\beta)b^{-1}),
$$
the product $c_1c_2^bc_3^{b^{-1}}$ is a commutator.%
\footnote{We are using the standard notation:
$x^y=y^{-1}xy$, $[x,y]=x^{-1}y^{-1}xy$.}
Hence
$$
c_1c_2c_3=c_1c_2^bc_3^{b^{-1}}
\Bigl((c_3^{-1})^{b^{-1}}c_2\Bigr)\Bigl((c_2^{-1})^bc_3\Bigr)
=(c_1c_2^bc_3^{b^{-1}})[c_2^{-1}c_3^{b^{-1}},b]
$$
is the product of two commutators.
\end{proof}



\bibliographystyle{amsplain}
\bibliography{/Volumes/Data/My_Math_Papers/bib}

\newcommand{\noopsort}[1]{} \newcommand{\singleletter}[1]{#1}
  \providecommand{\href}[2]{#2} \providecommand{\url}[1]{\texttt{#1}}
  \providecommand{\nolinkurl}[1]{\texttt{#1}}
\providecommand{\bysame}{\leavevmode\hbox to3em{\hrulefill}\thinspace}
\providecommand{\MR}{\relax\ifhmode\unskip\space\fi MR }
\providecommand{\MRhref}[2]{%
  \href{http://www.ams.org/mathscinet-getitem?mr=#1}{#2}
}
\providecommand{\href}[2]{#2}
\begin{thebibliography}{10}

\bibitem{AhlbrandtZiegler:1986:qfatct}
Gisela Ahlbrandt and Martin Ziegler, \emph{Quasi finitely axiomatizable totally
  categorical theories}, Ann.\ Pure Appl.\ Logic \textbf{30} (1986), no.~1,
  63--82.

\bibitem{BardakovTolstykh:2007:iaTgF}
Valery~G. Bardakov and Vladimir~A. Tolstykh, \emph{Interpreting the arithmetic
  in {T}hompson's group\/ ${F}$}, J.\ Pure Appl.\ Algebra \textbf{211} (2007),
  no.~3, 633--637, Preprint: \href{http://arxiv.org/abs/math/0701748}
  {\nolinkurl{arXiv:math/0701748}}.

\bibitem{BelkBrown:2005:fdeTgF}
James~M. Belk and Kenneth~S. Brown, \emph{Forest diagrams for elements of
  {T}hompson's group {$F$}}, Internat.\ J.\ Algebra Comput. \textbf{15} (2005),
  no.~5--6, 815--850.

\bibitem{BieriStrebel:pp1985:gPLhrl}
Robert Bieri and Ralph Strebel, \emph{On groups of ${PL}$-homeomorphisms of the
  real line}, Notes, Math.\ Sem.\ der Univ.\ Frankfurt, 1985.

\bibitem{BleakGGHMNS:pp2007:}
Collin Bleak, Alison Gordon, Garrett Graham, Jacob Hughes, Francesco Matucci,
  Hannah Newfield-Plunkett, and Eugenia Sapir, \emph{Using dynamics to analyze
  centralizers in the generalized {H}igman-{T}hompson groups\/ ${V}_n$},
  Incomplete preprint.

\bibitem{BrinSquier:2001:pcrcgplhrl}
Matthew~G. Brin and Craig~C. Squier, \emph{Presentations, conjugacy, roots, and
  centralizers in groups of piecewise linear homeomorphisms of the real line},
  Comm.\ Algebra \textbf{29} (2001), no.~10, 4557--4596.

\bibitem{Brown:1987:fpg}
Kenneth~S. Brown, \emph{Finiteness properties of groups}, J.\ Pure Appl.\
  Algebra \textbf{44} (1987), no.~1--3, 45--75.

\bibitem{CannonFloydParry:1996:inRTg}
James~W. Cannon, William~J. Floyd, and Walter~R. Parry, \emph{Introductory
  notes on {R}ichard {T}hompson's groups}, Enseign.\ Math.\ (2) \textbf{42}
  (1996), no.~3--4, 215--256, Preprint:
  \url{http://www.geom.uiuc.edu/docs/preprints/lib/GCG63/thompson.ps}.

\bibitem{Davis:1973:htpiu}
Martin Davis, \emph{Hilbert's tenth problem is unsolvable}, Amer.\ Math.\
  Monthly \textbf{80} (1973), 233--269.

\bibitem{DelonSimonetta:1998:uwpspsf}
Fran{\c c}oise Delon and Patrick Simonetta, \emph{Undecidable wreath products
  and skew power series fields}, J.\ Symbolic Logic \textbf{63} (1998), no.~1,
  237--246.

\bibitem{DennisVaserstein:1989:clg}
R.~Keith Dennis and Leonid~N. Vaserstein, \emph{Commutators in linear groups},
  $K$-Theory \textbf{2} (1989), no.~6, 761--767.

\bibitem{Ershov:1980:prkm-rus}
Yuri~L. Ershov, \emph{Problemy razreshimosti i konstruktivnye modeli
  [{D}ecision problems and constructivizable models]}, Matematicheskaya Logika
  i Osnovaniya Matematiki [{M}athematical Logic and Foundations of
  Mathematics], Nauka, Moscow, 1980 (Russian).

\bibitem{Godel:2006:fuSPMvS1-ger}
Kurt G{\"o}del, \emph{{\"U}ber formal unentscheidbare {S}{\"a}tze der
  {P}rincipia {M}athematica und verwandter {S}ysteme. {I} [{O}n formally
  undecidable propositions of {P}rincipia {M}athematica and related systems.
  {I}]}, Monatsh. Math. \textbf{149} (2006), no.~1, 1--30 (German), Reprinted
  from Monatsh.\ Math.\ Phys.\ {\bf 38} (1931), 173--198 [MR1549910], With an
  introduction by Sy-David Friedman. \MR{MR2260656}

\bibitem{Higman:1974:fpisg}
Graham Higman, \emph{Finitely presented infinite simple groups}, Notes on Pure
  Mathematics, vol.~8, Australian National University, Canberra, 1974.

\bibitem{Hodges:1993:mt}
Wilfrid Hodges, \emph{Model theory}, Encyclopedia of Mathematics and its
  Applications, vol.~42, Cambridge University Press, 1993.

\bibitem{Hodges:1997:smt}
\bysame, \emph{A shorter model theory}, Cambridge University Press, 1997.

\bibitem{Khelif:2007:bisQFAegrac-fr}
Anatole Kh{\'e}lif, \emph{Bi-interpr{\'e}tabilit{\'e}e et structures {QFA} :
  {\'e}tude de groupes r{\'e}solubles et des anneaux commutatifs
  [bi-interpretability and qfa structures: study of some soluble groups and
  commutative rings]}, C.\ R.\ Math.\ Acad.\ Sci.\ Paris \textbf{345} (2007),
  no.~2, 59--61. \MR{MR2343552 (2008e:03058)}

\bibitem{Matijasevich:1970:dpm-rus}
Yuri~V. Matijasevich, \emph{Diofantovost' perechislimykh mnozhestv [{T}he
  {D}iophantiness of enumerable sets]}, Dokl.\ Akad.\ Nauk SSSR \textbf{191}
  (1970), 279--282 (Russian), English translation in Soviet Math.\ Dokl.

\bibitem{MostowskiRobTar:1971:ueua}
Andrzej Mostowski, Raphael~M. Robinson, and Alfred Tarski, \emph{Undecidability
  and essential undecidability in arithmetic}, Undecidable theories, Studies in
  logic and the foundations of mathematics, North-Holland Publishing Co.,
  Amsterdam, 1971, pp.~36--74.

\bibitem{MostowskiTarski:1949:uaitr}
Andrzej Mostowski and Alfred Tarski, \emph{Undecidability in the arithmetic of
  integers and in the theory of rings}, J.\ Symbolic Logic \textbf{14} (1949),
  76.

\bibitem{Nies:2007:dg}
Andr{\'e} Nies, \emph{Describing groups}, Bull.\ Symbolic Logic \textbf{13}
  (2007), no.~3, 305--339.

\bibitem{Noskov:1983:etkpprg-rus}
Gennady~A. Noskov, \emph{Ob elementarnoj teorii konechno porozhdyennoj pochti
  razreshimoj gruppy [{O}n the elementary theory of a finitely generated almost
  solvable group]}, Izv.\ Akad.\ Nauk SSSR Ser.\ Mat. \textbf{47} (1983),
  no.~3, 498--517 (Russian), English translation in Math.\ USSR Izvestiya.

\bibitem{Noskov:1984:etfgasg-eng}
\bysame, \emph{On the elementary theory of a finitely generated almost solvable
  group}, Math.\ USSR Izvestiya \textbf{22} (1984), no.~3, 465--482, Translated
  from Russian.

\bibitem{Pillay:1983:ist}
Anand Pillay, \emph{An introduction to stability theory}, Oxford Logic Guides,
  vol.~8, The Clarendon Press, Oxford University Press, New York, 1983.

\bibitem{Poizat:1985:ctm-fr}
Bruno Poizat, \emph{Cours de th{\'e}orie des mod{\`e}les [{C}ourse on the
  theory of models]. {U}ne introduction {\`a} la logique math{\'e}matique
  contemporaine [{A}n introduction to contemporary mathematical logic]}, Nur
  al-Mantiq wal-Ma'rifah, Bruno Poizat, Lyon, 1985 (French).

\bibitem{Poizat:1987:gs-fr}
\bysame, \emph{Groupes stables [{S}table groups]. {U}ne tentative de
  conciliation entre la g{\'e}om{\'e}trie alg{\'e}brique et la logique
  math{\'e}matique [{A}n attempt at reconciling algebraic geometry and
  mathematical logic]}, Nur al-Mantiq wal-Ma'rifah [Light of Logic and
  Knowledge], vol.~2, Bruno Poizat, Lyon, 1987 (French).

\bibitem{Poizat:2000:cmt-eng}
\bysame, \emph{A course in model theory. {A}n introduction to contemporary
  mathematical logic}, Universitext, Springer-Verlag, 2000, Translated from
  French by Moses Klein and revised by the author.

\bibitem{Poizat:2001:sg-eng}
\bysame, \emph{Stable groups}, Mathematical Surveys and Monographs, vol.~87,
  American Mathematical Society, 2001, Translated from the 1987 French original
  by Moses Gabriel Klein.

\bibitem{Prest:1988:mtm}
Mike~Y. Prest, \emph{Model theory and modules}, London Mathematical Society
  Lecture Note Series, vol. 130, Cambridge University Press, 1988.

\bibitem{Robinson:1951:ur}
Raphael~M. Robinson, \emph{Undecidable rings}, Trans.\ Amer.\ Math.\ Soc.
  \textbf{70} (1951), 137--159.

\bibitem{Rothmaler:2000:imt-eng}
Philipp Rothmaler, \emph{Introduction to model theory}, Algebra, Logic and
  Applications, vol.~15, Gordon and Breach Science Publishers, Amsterdam, 2000,
  Translated and revised from the 1995 German original by the author.

\bibitem{Scanlon:2008:ifgfbN}
Thomas Scanlon, \emph{Infinite finitely generated fields are biinterpretable
  with\/ $\mathbb {N}$}, J.\ Amer.\ Math.\ Soc. \textbf{21} (2008), no.~3,
  893--908, Preprint: {\catcode`\~=11
  \url{http://math.berkeley.edu/~scanlon/papers/pc4sep07.pdf}}.

\bibitem{Sela:pp2006:dgg8s}
Zlil Sela, \emph{Diophantine geometry over groups {VIII}: stability}, Preprint,
  \href{http://arxiv.org/abs/math/0609096}{\nolinkurl{arXiv:math/0609096}},
  2006.

\bibitem{Stein:1992:gplh}
Melanie Stein, \emph{Groups of piecewise linear homeomorphisms}, Trans.\ Amer.\
  Math.\ Soc. \textbf{332} (1992), no.~2, 477--514.

\bibitem{Tarski:1971:gmpu}
Alfred Tarski, \emph{A general method in proofs of undecidability}, Undecidable
  theories, Studies in logic and the foundations of mathematics, North-Holland
  Publishing Co., Amsterdam, 1971, pp.~1--35.

\bibitem{Wagner:2000:sg}
Frank~O. Wagner, \emph{Stable groups}, Handbook of algebra, vol.~2,
  North-Holland Publishing Co., Amsterdam, 2000.

\bibitem{anonymous:2004:Tg40y}
\emph{Thompson's {G}roup at 40 {Y}ears}, 2004, Problem list of the workshop
  held 11--14 Jan., 2004, at American Institute of Mathematics, Palo Alto,
  California,
  \url{http://www.aimath.org/WWN/thompsonsgroup/thompsonsgroup.pdf}.

\end{thebibliography}


\end{document}